\magnification 1200
\input amssym.def
\input amssym.tex
\parindent = 40 pt
\parskip = 12 pt
\font \heading = cmbx10 at 12 true pt
 at 22 true pt
 at 19 true pt
 at 7 true pt
\def \R{{\bf R}}

\centerline{\heading Oscillatory Integral Decay, Sublevel Set Growth,}
\centerline{\heading and the Newton Polyhedron}
\rm
\line{}
\line{}
\centerline{\heading Michael Greenblatt}
\centerline{greenbla@uic.edu}
\line{}
\centerline{February 9, 2009}
\baselineskip = 12 pt
\font \heading = cmbx10 at 14 true pt
\line{}
\line{}
\noindent{\bf 1. Introduction}

\vfootnote{}{This research was supported in part by NSF grant DMS-0654073}In this paper we 
consider two types of integrals. Suppose $S(x)$ is a real-analytic function defined in a neighborhood
of the origin in $\R^n$. The first type of integral being considered are sublevel set integrals 
of the form
$$I_{S, \phi}(\epsilon) =  \int_{\{x : 0 < S(x) < \epsilon\}} \phi(x) \,dx \eqno (1.1a)$$
$$I_{|S|, \phi}(\epsilon) =  \int_{\{x : |S(x)| < \epsilon\}} \phi(x) \,dx \eqno (1.1b)$$
Here $\phi(x)$ is a smooth nonnegative
real-valued function supported within the domain of definition of $S(x)$ satisfying $\phi(0) > 0$.
Such integrals have been considered for example in [PSSt] and [Va], and are closely related to 
Gelfand-Leray functions.
We are interested in the behavior of $I_{S, \phi}(\epsilon)$ or $I_{|S|, \phi}(\epsilon) = 
I_{S, \phi}(\epsilon) + I_{-S, \phi}(\epsilon)$ as $\epsilon \rightarrow 0$.

\noindent The second type of integral under consideration are oscillatory integrals 
$$J_{S, \phi}(\lambda) = \int_{\R^n} e^{i \lambda S(x)} \phi(x)\,dx \eqno (1.2)$$
Again $\phi(x)$ is a smooth 
real-valued function supported within the domain of definition of $S(x)$, but we make no 
assumption of nonnegativity on $\phi(x)$.
Here we are interested in the behavior of $J_{S,\phi}(\lambda)$ as $|\lambda| \rightarrow \infty$. 
Since $\phi(x)$ is real, it suffices to consider the behavior of $J_{S,\phi}(\lambda)$ as 
$\lambda \rightarrow + \infty$. 

In this paper, extending the methods of [G1] we will prove theorems
generalizing a well-known theorem of Varchenko (Theorem 1.1 below) concerning oscillatory integrals 
$J_{S,\phi}$. They will be derived from analogous results proven here for the sublevel
integrals $I_{|S|,\phi}$. Varchenko's
theorem requires a certain nondegeneracy condition on the faces of the Newton polyhedron on $S$. In this paper,
we will show in Theorems 1.2 and 1.3 that the estimates he obtained also hold for a significant class of 
$S(x)$ for which this
nondegeneracy condition does not hold. Thus in problems where one wants to switch coordinates
to a coordinate system where Varchenko's estimates are valid, one has greater flexibility by
using the results of this paper. It should be pointed out that the methods of [G1] were influenced 
by those of [V] and therefore [V] can be viewed as an antecedent to this paper.

We will also exhibit some weaker estimates for more 
general situations, including some where
the estimates of Theorem 1.1 in fact do not hold. We will see that our conditions on $S(x)$ in 
Theorem 1.3 for Varchenko's estimates to hold are optimal
in some situations (Theorem 1.4). In two dimensions (Theorem 1.5), we will give
a characterization of the $S(x)$ for which the Newton polygon determines sharp estimates in
the fashion of Theorem 1.1; this too will hold for both the sublevel and oscillatory integrals.
This may be viewed as a generalization of [G3], at least for real-analytic phase. 

Integrals of the form $(1.1a)-(1.1b)$ and $(1.2)$ come up frequently
in analysis. For example, oscillatory integrals of the form $(1.2)$ arise in PDE's, mathematical
physics, and in harmonic analysis applications such as finding the decay of Fourier transforms of
surface-supported measures and associated problems concerning the restriction and Kakeya problems.
We refer to [AGV] chapter 6 and [S] chapter 8 for more information on such issues. The stability of 
oscillatory integrals of this kind under perturbations of the phase function $S(x)$ is related to a 
number of issues in complex geometry and has been
studied for example in [PSSt] and [V]. Also, operator versions of these oscillatory integrals have 
been extensively analyzed, for example in [G4] [G5] [GrSe] [R] [PS] [Se]. Furthermore, as will
be seen, our theorems concerning $I_{|S|, \phi}$ directly imply corresponding results for how the
measure of $\{x \in U : 0 <|S(x)| < \epsilon\}$ goes to zero as $\epsilon \rightarrow 0$. Here 
$U$ is a sufficiently small open set containing the origin. 
These come up for example in the analysis of Radon transforms such as in [C2] or [G5]. 

If $S(0) \neq 0$ and $\phi$ is supported on a sufficiently small neighborhood of the origin, then 
$I_{S, \phi}(\epsilon) = 0$ for small enough $\epsilon$  and thus is not interesting to analyze. 
In studying
$(1.2)$, one can always reduce to the case where $S(0) = 0$ by factoring out a $e^{i\lambda S(0)}$.
Hence it does no harm to assume that $S(0) = 0$ in the analysis of $J_{S,\phi}$ either. Furthermore,
if $\nabla S(0) \neq 0$, one easily
has that $I_{S, \phi}(\epsilon) \sim \epsilon$ as $\epsilon \rightarrow 0$ for $\phi$ supported near 
the origin. Also, by integrating by parts repeatedly in the $\nabla S(0)$ direction, one also has 
$|J_{S,\phi}| < C_N \lambda^{-N}$ as $\lambda \rightarrow +\infty$ if the support of $\phi$ is 
sufficiently small. Therefore the interesting situation for both $I_{S, \phi}$ and $J_{S, \phi}$ is when
$\nabla S(0) = 0$. Hence in this paper we will always assume that
$$S(0) = 0\,\,\,\,\,\,\,\,\,\, \nabla S(0) = 0 \eqno (1.3)$$
By Hironaka's resolution of singularities one has asymptotic expansions for both $I_{S, \phi}$ and 
$J_{S, \phi}$ if $\phi$ is supported in a sufficiently small neighborhood of the origin (see [G2] for 
elementary proofs). Namely, if $S(0) = 0$ and $\phi$ is supported in a sufficiently small
neighborhood of the origin one can asymptotically write
$$I_{S, \phi}(\epsilon) \sim \sum_{j = 0}^{\infty} \sum_{i=0}^{n-1} c_{ij}(\phi) \ln(\epsilon)^i \epsilon^
{r_j} \eqno (1.4a)$$
$$J_{S, \phi}(\lambda) \sim \sum_{j = 0}^{\infty} \sum_{i=0}^{n-1} d_{ij}(\phi) \ln(\lambda)^i \lambda^
{-s_j} \eqno (1.4b)$$
Here $\{r_j\}$ and $\{s_j\}$ are increasing arithmetic progressions of positive rational numbers 
independent of $\phi$ deriving from the resolution of singularities of $S$. Using
resolution of singularities one can show that the
smallest $r_j$ for which some $c_{ij}(\phi)$ is nonzero will not depend on what $\phi$ is, and similarly the 
largest $i$ for which $c_{ij}(\phi)$ is nonzero for this $j$ also is independent of $\phi$. 
(This uses the nonnegativity assumption on $\phi$ and that $\phi(0) > 0$). Hence as $\epsilon 
\rightarrow 0$, $I_{S, \phi}(\epsilon)$ will always be of the same order of magnitude.  Inspired by
terminology from the text [AGV], we refer to the value 
of $r_j$ in this case as the {\it growth index} of $S$ at the origin, and the 
corresponding value
of $i$ is referred to as the ${\it multiplicity}$ of this index. We define the growth index of $|S|$ to
be the minimum of the growth indices of $S$ and $-S$, with its multiplicity that of $S$ or
$-S$. The multiplicity taken to be the maximum of the
multiplicities of this growth index for $S$ and $-S$ if they both have the same growth index. Note that the above considerations 
imply that if $U$ is a sufficiently small neighborhood of the origin, then the measure of
$\{x \in U: 0 < S(x) < \epsilon\} \sim |\ln\epsilon|^i \epsilon^{r_j}$ as $\epsilon \rightarrow 0$, 
where $r_j$ is the growth index and $i$ is the
multiplicity of that index. As a result, knowing the growth index and its multiplicity gives the correct 
order of magnitude for such sublevel set volumes as $\epsilon \rightarrow 0$.

In the case of $J_{S, \phi}$, one
does not necessarily have that the smallest $r_j$ for which a $d_{ij}(\phi)$ is nonzero is the same for 
all $\phi$ (which is no longer even assumed to be nonnegative), so the above definition of index does 
not make sense. Instead, similar to [AGV] we define the {\it oscillation index} of $S$ at the origin 
to be the minimal 
$s_j$ for which for any sufficiently small neighborhood $U$ of the origin,
$d_{ij}(\phi)$ is nonzero for some $\phi$ supported in $U$. The multiplicity of this index $s_j$ 
is defined to be the maximal $i$ such that for any sufficiently small neighborhood $U$ of the origin
there is a $\phi$ supported on $U$ such that $d_{ij}(\phi)$ is nonzero for this minimal $s_j$.

In general, the growth or oscillation index and their multiplicities are determined by the zero set
of $S$ in a complicated
way. However, there are a number of situations when they can be determined from the Taylor series of
$S(x)$ at the origin in a nice geometric way, a fact discovered by Varchenko in [V]. Heuristically 
speaking, these situations correspond to when the zero of $S(x)$ at the origin is stronger than
any zero of $S(x)$ outside the coordinate hyperplanes $\{x_i = 0\}$. To indicate how the
index and its multiplicity are determined in these situations, we first define some terminology.

\noindent {\bf Definition 1.1.} Let $S(x) = \sum_{\alpha} s_{\alpha}x^{\alpha}$ denote the 
Taylor expansion of $S(x)$ at the origin.
For any $\alpha$ for which $s_{\alpha} \neq 0$, let $Q_{\alpha}$ be the octant $\{x \in \R^n: 
x_i \geq \alpha_i$ for all $i\}$. Then the {\it Newton polyhedron} $N(S)$ of $S(x)$ is defined to be 
the convex hull of all $Q_{\alpha}$.  

In general, a Newton polyhedron can contain faces of various dimensions in various configurations. 
These faces can be either compact or unbounded. 
In this paper as well as in [V], an important role is played by the following functions, defined for compact faces of
the Newton polyhedron. A vertex is always considered to be a compact face of dimension zero.

\noindent {\bf Definition 1.2.} Suppose $F$ is a compact face of the $N(S)$. Then
if $S(x) = \sum_{\alpha} s_{\alpha}x^{\alpha}$ denotes the Taylor expansion of $S$ like above, 
define $S_F(x) = \sum_{\alpha \in F} s_{\alpha}x^{\alpha}$

\noindent Also useful is the following terminology.

\noindent {\bf Definition 1.3.} Assume $S(x)$ is not identically zero. Then the 
{\it Newton distance} of $S(x)$ is defined to be $\inf \{t: (t,t,...,t,t) \in N(S)\}$.

The above-mentioned characterization in [V] of the oscillation index $S$ at $0$ and its 
multiplicity is as follows.

\noindent {\bf Theorem 1.1. (Varchenko)} Suppose for each compact face $F$ of $N(S)$, the function
$\nabla S_F (x)$ is nonvanishing on $(\R - \{0\})^n$. Further suppose that the Newton distance of 
$S$ is equal to some $d > 1$. Then the
oscillation index of $S$ at $0$ is given
by ${1 \over d}$. If the face of $N(S)$ (compact or not) that 
intersects the line $\{(t,t,...,t,t): t \in \R\}$ in its interior has dimension $k$, then the 
multiplicity of this index is given by $n - k - 1$.

For the purposes of Theorem 1.1, if the line $\{(t,t,...,t,t): t \in \R\}$ intersects $N(S)$ at a vertex,
then one takes $k = 0$. 

In this paper, we generalize Theorem 1.1 to a large class of functions where the $S_F(x)$ are not required
to have nonvanishing gradient, and prove analogues for the sublevel set integrals. We also prove weaker
substitutes for more degenerate situations including some when the conclusions of Theorem 1.1 do not 
necessarily hold. The methods of this paper are closely tied to the methods of [G1]. In turn, [G1] has
antecedents in the earlier two-dimensional algorithms [G4]-[G5], and also [PS] and [V]. There has
furthermore been much important work in sublevel set estimates and associated stability problems in the 
complex-analytic setting, such as in [DKo] [PSt1] [PSt2]. In [PSt1] and [PSt2], the method of algebraic 
estimates is used for this purpose; in [PSt2] resolution of singularities algorithms of Bierstone and 
Milman such as [BM] are also used. The complex methods tend to be rather different from the real ones since the
results obtainable in the complex case are quite a bit stronger than those obtainable in the real 
situation. 

In the theorems below, 
$S(x)$ is a real-analytic function, not identically zero, defined in 
a neighborhood of the origin and satisfying $(1.3)$. $d > 0$ denotes the Newton distance of $S(x)$. 
$C(S)$ denotes the face (compact or not) of $N(S)$ intersecting the line $\{(t,t,...,t,t): t 
\in \R\}$ in its interior, and $k$ denotes the dimension of $C(S)$. If the line intersects $N(S)$ 
at a vertex, we let $C(S)$ be this vertex and take $k = 0$. 

\noindent {\bf Theorem 1.2.} 

\noindent {\bf a)} As $\epsilon \rightarrow 0$, one has
$$I_{|S|,\phi}(\epsilon) >  C |\ln\epsilon|^{n - k - 1} \epsilon^{{1 \over d}}$$
\noindent {\bf b)} If for each compact face $F$ of $N(S)$ any zero of $S_F(x)$ in $(\R - \{0\})^n$ 
has order at most $d$, then as $\epsilon \rightarrow 0$ one has 
$$I_{|S|,\phi}(\epsilon) < C' |\ln\epsilon|^{n - k}\epsilon^{{1 \over d}}$$  
In this situation, as long as there is no compact face $F$ of $N(S)$ with $F \subset C(S)$ such that 
$S_F(x)$ has a zero 
of order $d$ somewhere in $(\R - \{0\})^n$, then one has the stronger estimate (compare with part {\bf a}) ) 
$$I_{|S|,\phi}(\epsilon) < C' |\ln\epsilon|^{n - k - 1}\epsilon^{{1 \over d}}$$
\noindent {\bf c)} If the maximum order of any zero of any $S_F(x)$ ($F$ compact) on 
$(\R - \{0\})^n$ is  $d' > d$, then as $\epsilon \rightarrow 0$ one can at least say that
$$I_{|S|,\phi}(\epsilon) < C' \epsilon^{{1 \over d'}}$$

We next come to our three-dimensional result. One can get somewhat stronger results
in three dimensions using a theorem of Karpushkin in [K] concerning the stability of growth indices
under deformations of the phase in $n - 1 = 2$ dimensions. A version of this theorem 
that sufficies for our purposes is as follows.

\noindent {\bf Theorem (Karpushkin)} Suppose $f(x_1,x_2,t_1,...,t_m)$ is a real-analytic function on a 
neighborhood of the origin in $\R^{m + 2}$ and $f(x_1,x_2,0,...,0)$ has growth index $c$ at the origin as 
a function of $x_1$ and $x_2$. Then for any $\mu > 0$, there is a constant $A_{\mu}$ and a 
neighborhood $U^{\mu}_1 \times U^{\mu}_2$ of the origin in the $(x,t)$ variables such that for 
$(t_1,...,t_m) \in U^{\mu}_2$ one has
$$|\{(x_1,x_2) \in U^{\mu}_1: |f(x_1,x_2,t_1,...,t_m)| < \epsilon\}| \leq A_{\mu}\epsilon^{c -
\mu}$$

To state our three-dimensional theorem, we need to
consider the growth index of a polynomial $S_F(x)$ at a point $a \neq 0$. By this we mean the
growth index of $S_F(x + a)$ at $x = 0$. When $S_F(a) \neq 0$, we define this growth index to infinity,
 and when $S_F(a) = 0$ but $\nabla S_F (a) \neq 0$, we take the growth index to be 1. 

\noindent {\bf Theorem 1.3.} Suppose $n = 3$. Then the following hold.

\noindent {\bf a)} As $\epsilon \rightarrow 0$, one has 
$$I_{|S|,\phi}(\epsilon) > C |\ln\epsilon|^{2 - k} \epsilon^{{1 \over d}}$$
\noindent {\bf b)} Suppose the growth index of every $|S_F|$ ($F$ compact) at any point in 
$(\R - \{0\})^3$ is at least ${1 \over d}$. Then as $\epsilon \rightarrow 0$ one has
$$I_{|S|,\phi}(\epsilon) < C' |\ln\epsilon|^2 \epsilon^{{1 \over d}}$$
If the growth index of every $|S_F|$ ($F$ compact) on $(\R - \{0\})^3$ is actually greater than ${1 \over 
d}$ at each point in $(\R - \{0\})^3$, then as $\epsilon \rightarrow 0$ one has the stronger (compare
with {\bf a}))
$$I_{|S|,\phi}(\epsilon) < C' |\ln\epsilon|^{2 - k}\epsilon^{{1 \over d}}$$
\noindent {\bf c)} Let $a$ denote the infimum over all compact faces $F$ of $N(S)$ and all $x \in 
(R - \{0\})^3$ of the growth index of $|S_F|$ at $x$. If $a  <  {1 \over d}$, then as $\epsilon
\rightarrow 0$ one has
$$I_{|S|,\phi}(\epsilon) < C' |\ln\epsilon|^2 \epsilon^{a}$$
We next have the following result, which may be viewed as a sort of converse to the type of result given in
Theorem 1.3, at least for the face $C(S)$. It holds in all dimensions.

\noindent {\bf Theorem 1.4.} Suppose $C(S)$ is a compact face of $N(S)$.

\noindent {\bf a)} Suppose there is some $x \in (R - \{0\})^n$ such that the growth index of 
$|S_{C(S)}|$ at $x$ is $a < {1 \over d}$. Then for some $a'  < {1 \over d}$, as $\epsilon
\rightarrow 0$ one has 
$$I_{|S|,\phi}(\epsilon) > C \epsilon^{a'}$$
\noindent {\bf b)} Suppose there is a $x \in (R - \{0\})^n$ such that $|S_{C(S)}|$ has a growth 
index of ${1 \over d}$ at $x$, with multiplicity $q$. Then as $\epsilon \rightarrow 0$ one has
$$I_{|S|,\phi}(\epsilon) > C |\ln\epsilon|^{n - k + q}\epsilon^{{1 \over d}}$$

In [V] it is shown that for any real-analytic phase in two dimensions, there are necessarily 
"adapted coordinates" in which the reciprocal of the Newton distance gives the correct oscillation index.
These results were generalized to smooth
phase in [IM]. There are many situations where the hypotheses of Theorems 1.2b do not
hold, but where they do hold after a coordinate change; take $S(x,y) = (x - y)^n$ in two-dimensions
for example. A natural question to ask is in which situations is there a 
coordinate change after which one is in the setting of Theorem 1.2b) or 1.3b). The two-dimensional 
proofs of [V] and [IM] use facts arising from two-dimensional resolution of singularities
such as Puiseux's theorem. Thus it would be reasonable to believe that proving analogues of such 
theorems in higher dimensions would use higher-dimensional resolution of singularities methods (and may 
be correspondingly more involved). 

In the other extreme, if one works in two dimensions and fixes a 
coordinate system, one has the following theorem, analogous to the results of [G3]. It will be 
a rather direct consequence of Theorems 1.2 and 1.4.

\noindent {\bf Theorem 1.5.} Suppose $n = 2$. Then the following hold.

\noindent {\bf a)} The growth index of $|S|$ at the origin is given by ${1 \over d}$ if and only if 
$C(S)$ is not a compact edge of $N(S)$ such that  
$S_{C(S)}$ has a zero on $(\R - \{0\})^2$ of order greater than $d$. If $C(S)$ is such a compact 
1-dimensional face, then the growth index is less than ${1 \over d}$. 

\noindent {\bf b)} When the growth index of $|S|$ at the origin is ${1 \over d}$, then the multiplicity
of this index is equal to $1 - k$, unless $S_{C(S)}$ has a zero on $(\R - \{0\})^2$ of order $d$,
in which case it is equal to $1$. 

By well-known methods relating sublevel integrals to oscillatory integrals, the above results about
the $I_{|S|,\phi}$ have direct implications for the $J_{S, \phi}$. Namely we have

\noindent {\bf Theorem 1.6.} 

\noindent {\bf a)} Suppose $\phi$ is nonnegative with $\phi(0) > 0$.

\noindent  If $d > 1$, or if $S(x)$ is either everywhere nonnegative or everywhere 
nonpositive in some neighborhood of the origin, then all statements and estimates analogous to 
those of Theorems 1.2-1.5 
hold for $J_{S,\phi}$ in place of $I_{|S|,\phi}$. If one is not in these situations, as long as the
growth index of $|S|$ is not an odd integer, then Theorems 1.2 and 1.3
hold for $J_{S,\phi}$ in place of $I_{|S|,\phi}$. In particular, they hold under any of the 
hypotheses of Theorem 1.2b) or 1.3b) if $d$ is not the reciprocal of an odd integer.

\noindent {\bf b)} For general smooth $\phi(x)$ and any $d$,
$J_{S,\phi}$ decays as fast or faster than the decay rates corresponding to any upper bound given by 
Theorems 1.2, 1.3 or 1.5 for $I_{|S|,\phi}$. 

\noindent {\bf Stability of Integrals.}

Karpushkin's theorem above can be described as a stability theorem for level set measures of 
two-dimensional integrals; he proved analogues for oscillatory integrals as well. The
analogues of these results in three or more dimensions do not hold, as exemplified by the following
result contained in [V].

\noindent {\bf Theorem [V].} Let $S_t(x,y,z) = (x^4 + tx^2 + y^2 + z^2)^2 + x^p + y^p + z^p$, where
$p \geq 9$. Then

\noindent {\bf a)} If $t > 0$, the oscillatory index of $S_t$ is ${3 \over 4}$. \hfill \break
\noindent {\bf b)} The oscillatory index of $S_0$ is ${5 \over 8}$.\hfill \break
\noindent {\bf c)} If $t < 0$, the oscillatory index of $S_t$ is given by ${1 \over 2} + \gamma(p)$,
where $\gamma(p) \rightarrow 0$ as $p \rightarrow \infty$.

The next two theorems are simple examples of this phenomenon that follow from Theorem 1.2; in particular
we avoid using the full Zariski three-dimensional resolution of singularities needed in [V] to prove the 
above result.   

\noindent {\bf Theorem 1.7.} Let $U_t(x,y,z) = x^4 + tx^2 + y^2 + z^2$.

\noindent {\bf a)} If $t > 0$, the growth index of $U_t$ at the origin is ${3 \over 2}$.\hfill \break
\noindent {\bf b)} The growth index  of $U_0$ at the origin is ${5 \over 4}$.\hfill \break
\noindent {\bf c)} If $t < 0$, the growth index $U_t$ at the origin is 1.

\noindent Theorem 1.7 will quickly lead to the following oscillatory integral analogue.

\noindent {\bf Theorem 1.8.} Let $V_t(x,y,z) = (x^4 + tx^2 + y^2 + z^2)^2$

\noindent {\bf a)} If $t > 0$, the oscillatory index of $V_t$ is ${3 \over 4}$.\hfill \break
\noindent {\bf b)} The oscillatory index of $V_0$ is ${5 \over 8}$.\hfill \break
\noindent {\bf c)} If $t < 0$, the oscillatory index of $V_t$ is ${1 \over 2}$.

\noindent {\bf Proofs of Theorem 1.7 and 1.8.}

If $t \geq 0$, then for each compact 
face $F$ of $N(U_t)$ the corresponding polynomial $(U_t)_{F}(x,y,z)$ has no zeroes on $(\R - \{0\})^3$. Hence
the growth index of $U_t$ in these situations is given by Theorem 1.2a)-b). Computing the
Newton distances, one sees that the growth index at the origin is equal to ${3 \over 2}$ if $t > 0$ and equal to 
${5 \over 4}$ when $t = 0$. This gives parts a) and b) of Theorem 1.7. Now assume $t < 0$. Since each
polynomial $(U_t)_{F}(x,y,z)$ has zeroes of order at most 1 on $(\R - \{0\})^3$, Theorem 1.2c)
implies that the growth index of $U_t$ is at least 1. To show it is exactly 1, we do a 
variable change, writing $(x,y,z) = (x,xy',xz')$. In the new coordinates $U_t(x,y,z)$ becomes the
function $W_t(x,y',z') = x^2(x^2 + t + (y')^2 + (z')^2)$. This has zeroes on the sphere $x^2 + (y')^2 +
(z')^2 = -t$, and $\nabla W_t(x,y',z')$ is nonzero at any such zero with $0 < |x| < \sqrt{-{t \over 2}}$. Going 
back into the $(x,y,z)$ coordinates, this means $U_t(x,y,z)$ has zeroes arbitrarily close to the 
origin at which
$\nabla U_t(x,y,z) \neq 0$. In a small neighborhood of each such zero, the measure of $\{(x,y,z): 
U_t(x,y,z) < \epsilon\}$ is bounded below by $C \epsilon$. Hence the growth index of $U_t$ is at 
most 1. We conclude that the growth index of $U_t$ is exactly 1, giving part c) of Theorem 1.7 and 
completing the proof of that theorem. 

We move to Theorem 1.8. The 
growth index at the origin of $V_t$ is half that of $U_t$. So if $t > 0$, the growth index is ${3 
\over 4}$, if $t = 0$ it is ${5 \over 8}$, and if $t < 0$ it is ${1 \over 2}$. We will see in  
the last paragraph of section 5 that if the growth index is less than 1 then the oscillatory index and 
the growth index are the same (this also follows pretty directly from Ch 7 of [AGV]). The desired
properties immediately follow and we are done.

\noindent {\bf 2. Geometric constructions from the Newton polyhedron} 

In this section we do a number of geometric constructions which will used
in later sections in proving the various estimates
of this paper. As indicated above, they are based on the resolution of singularities methods of [G1].
However, we do not need a full-fledged resolution of singularities algorithm for the purposes of this
paper. 

Heuristically speaking, what we will do is as follows. Suppose $S(x)$ is a real-analytic function
defined on a neighborhood of the origin. We will take a small neighborhood of the origin, and divide it 
(modulo sets of measure zero) into open slivers $W_{ij}$ whose closures each contains the origin. 
Each $W_{ij}$ 
corresponds to one vertex or compact face $F_{ij}$ of $N(S)$ of dimension $i$ in the sense that on $W_{ij}$,
the monomials $x^v$ for $v$ a vertex of $N(S)$ on $F_{ij}$ dominate the monomials $x^v$ for $v \in N(S)$ 
not on $F_{ij}$. Lemmas 2.0 and 2.1 make this notion precise. 

Next, each $W_{ij}$ will be further subdivided, modulo sets of measure zero, into open slivers $W_{ijp}$ 
to each of which there will be assigned an invertible map $\beta_{ijp}:Z_{ijp} \rightarrow W_{ijp}$. Each 
component
function of each $\beta_{ijp}$ is plus or minus a monomial in $x_1^{{1 \over N}}... x_n^{{1 \over N}}$
for some integer $N$, and for $i > 0$ each domain $Z_{ijp}$ satisfies inclusions of the form 
$$(0,\eta ')^{n-i} \times D_{ij} \subset Z_{ijp} \subset (0,1)^{n-i} \times D_{ij} \eqno (2.0)$$
Here $D_{ij}$ is a bounded open set whose closure is contained in $\{(x_{n-i+1},...,x_n): x_k > 0$ for
all $k \}$. Furthermore, the map $\beta_{ijp}$ is such that $S_{F_{ij}} \circ \beta_{ijp}(x)$ can be expressed as
$m(x_1,...,x_{n-i})T(x_{n-i+1},...,x_n)$, where $m(x_1,...,x_{n-i})$ is a monomial in the first $n-i$ variables. As
a result, a condition that the zeroes of $S_{F_{ij}}$ on $(\R - \{0\})^n$ are of order less than $d$
implies that the same condition holds for $T(x_{n-i+1},...,x_n)$. Similarly, the various 
other conditions stipulated on $S_{F_{ij}}(x)$ in the different lemmas imply that the same
condition holds for $T(x_{n-i+1},...,x_n)$. In addition, since the 
$x_{n-i+1},...,x_n$ variables are bounded above on $Z_{ijp}$ (by the boundedness of the $D_{ij}$),
the function $S_{F_{ij}} \circ \beta_{ijp}(x)$ is bounded above by $C m(x_1,...,x_{n-i})$. 
Analogously, using that the points in $D_{ij}$ have coordinates bounded below away from zero, a local
nonvanishing $p$th derivative condition on $S_{F_{ij}}(x)$ will imply the corresponding $p$th
derivative of $T(x_{n-i+1},...,x_n)$ is bounded below on some open set, so that this derivative of 
$S_{F_{ij}} \circ \beta_{ijp}(x)$ is bounded below by $C' m(x_1,...,x_{n-i})$ in some
neighborhood. These facts are proven via the constructions of Theorem 2.2 and Lemmas 2.3 and 2.4. 

Because the terms $x^v$ for $v \in F_{ij}$ dominate on $W_{ijp}$, the difference $|S_{F_{ij}} \circ 
\beta_{ijp}(x) - S \circ \beta_{ijp}(x)|$ is bounded above by $\epsilon\, m(x_1,...,x_{n-i})
$ for a small $\epsilon$ (Lemma 2.1). So one has $|S\circ \beta_{ijp}(x)| < C''m(x_1,...
x_{n-i})$. One can also do the constructions are done so that any given derivative of 
$S \circ \beta_{ijp}(x)$ is a also a small perturbation of the corresponding derivative of
$S_{F_{ij}} \circ \beta_{ijp}(x)$; if $S_{F_{ij}} \circ \beta_{ijp}(x)$ satisfies a nonvanishing $p$th 
derivative condition on a small open set, then so does $S \circ \beta_{ijp}(x)$ which is therefore
bounded below by $C'''m(x_1,...,x_{n-i})$. This enables one to use van der 
Corput type lemmas in the $x_{n-i+1}....x_n$ variables to prove various desired estimates. 
The most convenient such van der Corput lemma for our purposes is that of [C1], which says that if $f(x)$ 
is a $C^{k+1}$ function on an interval $I$ whose $k$th derivative is bounded below by $\eta$, then one 
has 
$$|\{x \in I: |f(x)| < \delta\}| < C_k ({\delta \over \eta})^{1 \over k} \eqno (2.1)$$

Note that the properties being used
here are quite a bit weaker than those of a full resolution of singularities theorem since we only
have an upper bound for the blown-up function and a lower bound for its derivative in terms of a 
monomial (the blown-up function can even have a complicated zero set), but this suffices for our purposes.

For the three dimensional result, instead of getting uniform estimates from a Van der Corput-type lemma,
one considers the growth index directly. One uses Karpushkin's theorem to show that locally the growth 
index of $S \circ \beta_{ijp}(x)$ is the same as that of $S_{F_{ij}} \circ \beta_{ijp}(x)$, which in turn
is the same as that of  $T(x_3)$ or $T(x_2,x_3)$ (for $i = 1$ and 2 respectively.) One always has has 
to be careful that perturbing $S_{F_{ij}} \circ \beta_{ijp}(x)$ into $S \circ \beta_{ijp}(x)$ can be 
done in such a way that Karpushkin's result applies. 

To enable us to use van der Corput lemmas most effectively, one should have a good idea of what the monomials 
$m(x_1,...,x_{n-i})$ are. Fortunately, Lemmas 2.5 and 2.6 give us a way of doing this. Namely, if one 
redefines $\beta_{ijp}(x)$ such that for $q \leq n - i$ one replaces each $x_q$ by $x_q^{l_q}$, where the 
$l_q$ are chosen such that the determinant of $\beta_{ijp}(x)$ is constant, then each
variable $x_q$ for $q \leq n-i$ appears to at most the $d$th power in $m(x_1,...,x_{n-i})$, where as 
usual $d$ is the Newton distance. Furthermore, the $d$th power appears at in most 
$n - k$ variables, where $k$ is the dimension of the central face $C(S)$, and it appears $n - k$ times 
if and only if $F_{ij} \subset C(S)$. These things are proven in Lemma 2.6. One then proves the 
estimates of Theorems 1.2b-c by first using the appropriate Van der Corput-type lemma in 
a direction in the $x_{n-i+1},...,x_n$ variables, then taking absolute values and integrating in the 
remaining $x_{n-i+1},...,x_n$ variables, and then integrating the resulting function of the first $n - i$
variables. As one might guess, one needs to take a lot of care in carrying out this strategy. 

It should be pointed out that in the above description, we always assumed $i > 0$. But there are also
$W_{0jp}$; fortunately these are easy to deal with since the functions $T(x_{n-i+1},...,x_n)$ are 
replaced by a constant. 

To give a concrete and easy-to-understand example of the above considerations, in three dimensions consider the function
$S(x,y,z) = x^2 + y^2 - z^2$. Then the Newton distance of $S$ is ${2 \over 3}$, and the 
functions $S_e(x,y,z)$ either have no zeroes on $(\R - \{0\})^3$, or have zeroes of order 1 on 
$(\R - \{0\})^3$. In the above language, this says that the exponents
appearing in each monomial $m(x_1)$ or $m(x_1,x_2)$ are at most ${2 \over 3}$, while the functions 
$T(x_2,x_3)$ or $T(x_3)$ can have zeroes of order as high as 1. We focus our attention on the 
situation where $F_{ij}$ is the main 2-dimensional face; the 1-dimensional faces where $S_e(x,y,z)$ 
has a zero will behave similarly to the following. Then $i = 2$, and $T(x_2,x_3)$ has zeroes of order 1. 
By first using the Van der Corput lemma $(2.1)$ in an appropriate direction in the $x_2 x_3$ variables, 
then integrating in the orthogonal $x_2 x_3$ direction, and lastly integrating in $x_1$, using 
$(2.0)$ one gets that for some positive $\delta$ and $\delta'$, $|\{(x_1,x_2,x_3): |S \circ
\beta_{ijp}(x_1,x_2,x_3)|  < \epsilon\}|$ is comparable to $\int_0^{\delta} \max(\delta',{\epsilon \over {x_1^{2 
\over 3}}})\,dx_1 \sim \epsilon$. These are the weaker bounds of Theorem 1.2c). It is worth pointing out that
the oscillation is index here is the value ${3 \over 2}$ given by the Newton polyhedron since the
phase has nonvanishing Hessian. This is an example where one gets a smaller growth index 
(which is in fact 1 in this example) than oscillation index; by Theorem 1.6 for this
to happen $d$ must be less than 1. 

Next, suppose that instead of $S(x,y,z) = x^2 + y^2 - z^2$, one chooses $S(x,y,z) = x^4 + y^4 - z^4$. 
Then the Newton distance doubles to ${4 \over 3}$, yet the maximum order of any zero of any 
$S_e(x,y,z)$ is still 1. Since $1 < {4 \over 3}$, the stronger results of Theorem 1.2b) apply (It doesn't
help in this particular situation to use Theorem 1.3). Instead of ending out with an integration of
$\int_0^{\delta} \max(\delta',{\epsilon \over {x_1^{2 \over 3}}})\,dx_1$, one ends out with an integration of 
$\int_0^{\delta} \max(\delta', {\epsilon \over {x_1^{4 \over 3}}})\,dx_1$. Since the exponent in the
denominator is now greater than 1, the result is now comparable to $\epsilon^{3 \over 4}$. Simply put,
the zero of $S(x,y,z)$ at the origin now dominates the zeroes of $S(x,y,z)$ away from the origin on
$W_{ijp}$, so the Newton polyhedron now determines the growth index. On the other hand, in the previous 
example the reverse was true, so that the zeroes of $T(x_2,x_3)$ and its analogues from the other $W_{ijp}$
force the growth index to be smaller. Theorem 1.2c) says that, like in this example, that the growth index 
is bounded below by the reciprocal of the maximal order of a zero of the functions $T(x_2,x_3)$ or 
$T(x_3)$. In general, when the Newton polyhedron determines the growth index, the powers of at least 
one variable appearing in the integration for one $Z_{ijp}$ will have exponent at least 1, while when
the zeroes are too strong for that, all powers of all variables will be less than one. So these two 
examples, however simple, are fairly indicative. 

The lower bounds of Theorems 1.2a) and 1.3a) are not affected by the behavior of the zeroes of the various
$S_{F_{ij}}(x_1,...,x_n)$ since the zeroes can only cause one to obtain worse estimates than those given by
the Newton polyhedron. Thus in proving the lower bounds one can just restrict attention to some small 
subregion of $D_{ij}$ away from the
zeroes of the associated $T(x_{n-i+1},...,x_n)$. On this region $S \circ \beta_{ijp}(x) \sim 
m(x_1,...,x_{n-i})$ and the lower bounds determined by the Newton polyhedron are readily proven. 

\noindent We now begin proving our various lemmas. 

\noindent {\bf Lemma 2.0.} (Lemma 3.2 of [G1])
Let $v(S)$ denote the set of vertices of $N(S)$. There are $A_1, A_2 > 1$ such that if $C_0,...,C_n$
are constants with $C_0 > A_1$ and $C_{i+1} > C_i^{A_2}$ for all $i$, then one can define the $W_{ij}$ 
so that

\noindent {\bf a)} Let $i < n$. If the following two statements hold, then  $x \in W_{ij}$.

\noindent {\bf 1)} If $v \in v(S) \cap F_{ij}$ and  $v' \in v(S) \cap
(F_{ij})^c$ we have $ x^{v'} < C_n^{-1} x^v$. \parskip = 0pt

\noindent {\bf 2)} For all $v, w \in v(S) \cap F_{ij}$ we have $C_i^{-1}x^w < x^v < C_ix^w$. \parskip = 12pt

\noindent {\bf b)} There is a $\delta > 0$ depending on $N(S)$, and not on $A_1$ or $A_2$, such that if 
$x \in W_{ij}$, then the following two statements hold.

\noindent {\bf 1)} If $v \in v(S) \cap F_{ij}$ and  $v' \in v(S) \cap
(F_{ij})^c$ we have $ x^{v'} < C_{i+1}^{-\delta} x^v$. \parskip = 0pt

\noindent {\bf 2)} For all $v, w \in v(S) \cap F_{ij}$ we have $C_i^{-1}x^w < x^v < C_ix^w$. \parskip = 12pt

Informally, this gives a way of saying that the vertices of $F_{ij}$ dominate the Taylor series
of $S$ when $x \in W_{ij}$. Another way of making this precise is the following lemma.

\noindent {\bf Lemma 2.1.} Suppose $x \in W_{ij}$. Let $V \in v(S)$ be such that $x^V \geq x^v$ for
all $v \in v(S)$; if there is more than one such vertex let $V$ be any of them. Then if
$A_1$ is sufficiently large and $\eta$ is sufficiently small, for any positive $d$ one has the following estimate:
$$\sum_{\alpha \notin F_{ij}} |s_{\alpha}| |\alpha|^d x^{\alpha} < K(C_{i+1})^{-\delta''}x^V$$
Here $K$ is a constant depending on $d$ as well as the function $S(x)$, and $\delta'' > 0$ is a constant
depending on the Newton polyhedron of $S$. 

\noindent {\bf Proof.} There are several one-dimensional faces of $N(S)$ that contain $V$, and there are 
vectors $w_1,...w_N$ so that a given edge is given by $V + tw_l$ for a set of
nonnegative $t$. If any component of a vector $w_l$ is negative, the corresponding edge will terminate
at a vertex which we denote by $v_l$. Rescaling $w_l$ if necessary, we can assume that $v_l = V + w_l$.
If all components of a $w_l$ are nonnegative, then the edge is an
infinite ray. (It is not hard to show that $w_l$ is in fact some unit coordinate vector ${\bf e}_m$). In 
this situation we define $v_l = V + w_l$. Consequently, for all $l$ we have
$$x^{v_l} = x^V x^{w_l} \eqno (2.2)$$
I claim that, shrinking $\eta$ if necessary, we may assume that for all $l$ such that $v_l \notin F_{ij}$ we have
$$x^{v_l} < (C_{i+1})^{-\delta}x^V \eqno (2.3)$$
This is true if $v_l$ is a vertex of $N(S)$ by Lemma 2.0 above. It is true if $v_l$ is not a 
vertex since ${x^{v_l} \over x^V} = x^{w_l}$, which can be made less than $(C_{i+1})^{-\delta}$ by
shrinking $\eta$ appropriately since $w_l$ has only nonnegative components. So we can assume $(2.3)$ 
holds. Next, note that since $N(S)$ is a convex polyhedron we have
$$N(S) \subset \{V + \sum_{l=1}^N t_l w_l: t_l \geq 0\} \eqno (2.4)$$
For a positive integer $k$, define $B_k$ to be the set of points $\alpha$ with 
integer coordinates that are in $N(S)$ but not on $F_{ij}$ such that $\alpha$ can be written as $V +
\sum_l t_l w_l, t_l \geq 0$ with $k - 1 < \sum_{v_l \notin F_{ij}} |t_l| \leq k$. 
Let $E$ be a separating 
hyperplane for $N(S)$ such that $E \cap N(S) = F_{ij}$. Since $F_{ij}$ is bounded, we may let $a$
be a vector normal to $F_{ij}$ such that each component of $a$ is positive. For each $w_l$ not parallel to 
$F_{ij}$, the vector $w_l$ points "inward"; that is, $a \cdot w_l > 0$. Consequently, for a constant
$C$ depending only on $N(S)$, the points in $B_k$ are contained in the points of $(\R^+)^n$ between
$E$ and its translate $E + Ck a$. In particular each coordinate of a point in $B_k$ is bounded by 
$Ck$ and there at most $Ck^n$ of them. Next, writing a 
given $\alpha \in B_k$ as $V + \sum_l t_l w_l$ with $k - 1< \sum_{v_l \notin F_{ij}} |t_l| <  k$, we 
have
$$x^{\alpha} = x^V \prod_{l}(x^{w_l})^{t_l} = x^V \prod_{v_l \in F_{ij}}(x^{w_l})^{t_l}
\prod_{v_l \notin F_{ij}}(x^{w_l})^{t_l} \leq x^V \prod_{v_l \notin F_{ij}}(x^{w_l})^{t_l} \eqno (2.5)$$
The last inequality follows from $(2.2)$ and the maximality of $x^V$. Using $(2.3)$ and the definition
of $B_k$ we have
$$x^V \prod_{v_l \notin F_{ij}}(x^{w_l})^{t_l} <  x^V(C_{i+1})^{-\delta \sum_{v_l \notin F_{ij}} t_l}
< (C_{i+1})^{-\delta(k - 1)} x^V\eqno (2.6)$$
When $k = 1$, one has an inequality
$$x^{\alpha} < x^V (C_{i+1})^{-\delta \sum_{v_l \notin F_{ij}} t_l} <  (C_{i+1})^{-\delta '} x^V
\eqno (2.7)$$
Here $\delta '$ is the minimum of the finitely many positive numbers $\delta \sum_{v_l \notin 
F_{ij}} t_l$ that can appear in the right hand side of $(2.7)$. Since $S$ is real analytic, the 
coefficients $s_{\alpha}$ satisfy $|s_{\alpha}| < CM^{|\alpha|}$ for some $M$. Since the components
of any $\alpha$ in any $B_k$ are at most $Ck$, we have
$$ |s_{\alpha}| < C'M^k \eqno (2.8)$$
Since there are most $Ck^n$ points with integer coordinates in any $B_k$, inserting $(2.8)$ in $(2.6)$
or $(2.7)$ and adding gives the following for $k > 1$.
$$\sum_{\alpha \in B_k} |s_{\alpha}| |\alpha|^d x^{\alpha} < C' k^{d+n}M^{k} C_{i+1}^{-\delta(k - 1)} x^V$$
$$= C'M^2k^{d+n}C_{i+1}^{-\delta}(MC_{i+1}^{-\delta})^{(k-2)}x^V\eqno (2.9a)$$
If $k = 1$ we have
$$\sum_{\alpha \in B_1} |s_{\alpha}| |\alpha|^d x^{\alpha} < C' M C_{i+1}^{-\delta '}x^V \eqno (2.9b)$$
Adding this over all $k$, as long as $A_1^{\delta} > 2M$ so that each $MC_{i+1}^{-\delta} < {1 \over 2}$, we get
$$\sum_{\alpha \notin F_{ij}} |s_{\alpha}| |\alpha|^d x^{\alpha} < C'' C_{i+1}^{-\delta ''}x^V
\eqno (2.10)$$
Here $\delta'' = \min(\delta,\delta')$. This gives the lemma and we are done.

\noindent {\bf Corollary.} There is a constant $C$ such that on a sufficiently small neighborhood of 
the origin $|S(x)| \leq C \sum_{v \in v(S)} x^v$.

\noindent {\bf Proof.} It suffices to prove the corollary on a given $W_{ij}$. We have
$$|S(x)| \leq \sum_{\alpha} |s_{\alpha}| x^{\alpha} = \sum_{\alpha \in F_{ij}} |s_{\alpha}|
x^{\alpha} + \sum_{\alpha \notin F_{ij}} |s_{\alpha}| x^{\alpha}< C_0x^V + K(C_{i+1})^{-\delta '}x^V$$
$$ = (C_0 + K(C_{i+1})^{-\delta '})x^V \leq (C_0 + K(C_{i+1})^{-\delta '})\sum_{v \in v(S)} x^v$$
The corollary follows.

For the purposes of this paper, we need to do a further subdivision of a given $W_{ij}$ into finitely
many pieces $W_{ijp}$. The relevant properties of the $W_{ijp}$ are encapsulated by the following 
theorem.

\noindent {\bf Theorem 2.2.} If $A_1$ and $A_2$ are sufficiently large, each $W_{ij}$ can be, modulo
a set of measure zero, written as the union of finitely many open nonempty sets $W_{ijp}$ to each of
which is associated a bijective map $\beta_{ijp}: Z_{ijp}  \rightarrow W_{ijp}$ depending
on $N(S)$ and $(i,j,p)$, but not the particular subdivision being done, such that each component of 
$\beta_{ijp}(z)$ is a monomial in $(z_1^{{1 \over N}},...,z_n^{{1 \over N}})$ for some $N$, and such 
that for some $\mu ' > 0$ that is allowed to depend on the particular subdivision we have 

\noindent {\bf a)} When  $i = 0$, $(0, \mu ')^n \subset Z_{ijp} \subset (0,1)^n$.\parskip = 3 pt

\noindent {\bf b)} When $i > 0$, there are sets $D_{ij} \subset (C_i^{-e}, C_i^e)^i$ for some $e > 0$ 
depending on $N(S)$ such that $(0, \mu ')^{n-i} \times D_{ij} \subset Z_{ijp} \subset 
(0, 1)^{n-i} \times D_{ij}$

\noindent {\bf c)} When $i > 0$, write $z \in \R^n$ as $(\sigma,t)$ where
$\sigma \in \R^{n-i}$ and $t \in \R^i$. For any $v \in N(S)$, denote by $\sigma^{v'}t^{v''}$ the
function in $z$ coordinates that $x^v$ transforms into under the $x$ to $z$ coordinate change. When
$i = 0$, write $z = \sigma$ and for $v \in N(S)$ denote by $\sigma^{v'}$ the the function
$x^v$ transforms into. 
Then for any $v_1,v_2 \in F_{ij}$ we have $v_1' = v_2'$, while if $v_1 \in F_{ij}$
and $v_2$ is in $N(S)$ but not in $F_{ij}$, then $(v_2')_k \geq (v_1')_k$ for all $k$ with at least 
one component strictly greater. \parskip = 12 pt

The proof of Theorem 2.2 is very similar to the arguments of section 4 of [G1]. However, there 
are enough differences that we prove it separately here. We will do it through some constructions 
resembling Lemmas 4.1-4.3 of [G1], after which we will prove Theorem 2.2.

For each $i$ and $j$ let $f_{ij}$ be any vertex on on $F_{ij}$.
Since the face $F_{ij}$ is of dimension $i$, we may let $\{P_l\}_{l=1}^{n-i}$
be separating hyperplanes for $N(S)$ such that $F_{ij} = \cap_{l=1}^{n-i} P_l$. We write these 
hyperplanes as 
$$P_l = \{x: a^l \cdot  x = c^l\}$$
We can assume the $a^l$ have rational coefficients. The hyperplanes satisfy
$$N(S) \subset \cap_{l=1}^{n-i} \{x: a^l \cdot x \geq c^l\} \eqno (2.11)$$
Since $\cap_{m=1}^n\{x: x_m  \geq f_{ijm}\} \subset N(S)$, we also have
$$\cap_{m=1}^n\{x: x_m  \geq f_{ijm}\} \subset \cap_{l=1}^{n-i} \{x: a^l \cdot x \geq c^l\} \eqno 
(2.12)$$
Since $a^l \cdot f_{ij} = c^l$ for all $l$, if we shift $x$ in $(2.11)$ by $-f_{ij}$ we get
$$\cap_{m=1}^n\{x: x_m  \geq 0 \} \subset \cap_{l=1}^{n-i} \{x: a^l \cdot x \geq 0\} \eqno 
(2.13)$$
In the case where $i > 0$, we would like to extend the hyperplanes $a^l \cdot x = 0$ to a collection of 
$n$ independent hyperplanes such that 
$$ \cap_{m=1}^n\{x: x_m  \geq 0\} \subset \cap_{l=1}^n \{x: a^l \cdot x \geq 0\} 
\eqno (2.14)$$
(Note that $(2.14)$ is $(2.13)$ when $i = 0$.)
We do this by defining 
$a^l$ for $i < l < n$ to be unit coordinate vectors such that $a^1,...,a^n$ are linearly
independent. Once we do this, we have 
$$ \cap_{m=1}^n\{x: x_m  \geq 0\} \subset \cap_{l=n-i+1}^{n}\{x: a^l \cdot x \geq 0\}\eqno (2.15)$$
Combining with $(2.13)$ gives $(2.14)$.

Since the $a^l \cdot x \geq 0$ are $n$ independent hyperplanes intersecting at the origin, any $n-1$ 
of the hyperplanes 
intersect along a line through the origin. Write the directions of these lines as $b_l$, chosen so
that the $b_l$ have rational components and $a_l \cdot b_l > 0$. The $b_l$ span $\R^n$, so we may write
the $m$th unit coordinate vector ${\bf e}_m$ in the form
$${\bf e}_m = \sum_{l=1}^n d_{lm} b_l \eqno (2.16)$$
\noindent {\bf Lemma 2.3.} The coefficients $d_{lm}$ are all nonnegative rational numbers.

\noindent {\bf Proof.} 
By definition of $b_l$, we have
$$ \cap_{l=1}^n \{x: a^l \cdot x \geq 0\} = \{s: s = \sum_{p=1}^n s_p b_p \hbox { with } s_p \geq 0 \} 
\eqno (2.17)$$ 
Since each ${\bf e}_m$ is in $\cap_{m=1}^n\{x: x_m  \geq 0\} \subset \cap_{l=1}^n \{x: a^l \cdot x 
\geq 0\}$, $(2.17)$ says that each $d_{lm}$ 
is nonnegative. Elementary linear algebra gives a formula for the $d_{lm}$ which shows that they are
rational. This completes the proof.

We now do a coordinate change on each $W_{ij}$ for $i > 0$. Denoting the original coordinates of a point $x$ 
by $(x_1,...,x_n)$, we let the new coordinates be denoted by $(y_1,...,y_n)$, where 
$$y_m = \prod_{l=1}^n x_l^{d_{lm}} \eqno (2.18)$$
Observe that a monomial $x^{\alpha}$ becomes $y^{L(\alpha)}$ in the new coordinates, where $L$
is the linear map such that $L(b_l) = {\bf e}_l$ for all $l$. If $\bar{f}_{ij} = 
(\bar{f}_{ij1},...,\bar{f}_{ijn})$ denotes $L(f_{ij})$, then each $\bar{f}_{ijk} \geq 0$ since
each $d_{lm}$ is nonnegative. Furthermore, $L$ takes each hyperplane $P_l$ to
$\{y: y_l = \bar{f}_{ijl}\}$. Notice that each point $p$ of $F_{ij}$ is on $P_l$ for $l \leq n - i$. 
This means that the $l$th component of $L(p)$ is equal to $\bar{f}_{ijl}$ for $l \leq n - i$. 
So if $v$ and $v'$ are vertices of $N(S)$ on $F_{ij}$, the first $n-i$ components of $L(v - v')$ 
are zero. Hence  ${y^{L(v - v')}}$ is a function of
the last $i$ $y$-variables only. Write $y = (s,t)$, where $s$ is the first $n-i$ variables and $t$ 
is the last $i$ variables. Similarly, write $L= (L_1,L_2)$, where $L_1$ is the first $n-i$ 
components and $L_2$ is 
the last $i$ components. Recall from Lemma 2.0 that for any such $v$ and $v'$, any  
$x \in W_{ij}$ satisfies the inequalities
$$C_i^{-1} < {x^{v - v'}} < C_i \eqno (2.19a)$$
In terms of the $t$ variables this translates as 
$$C_i^{-1} < {t^{L_2(v - v')}} < C_i \eqno (2.19b)$$
Write $\log(t) = (\log(t_1),\log(t_2),...,\log(t_n))$. Equation $(2.19b)$ becomes
$$-\log(C_i) < \log(t) \cdot L_2(v - v') < \log(C_i) \eqno (2.20)$$
Since the set of all possible $L_2(v - v')$ for $v$ and $v'$ vertices of $S$ on $F_{ij}$ spans an 
$i$-dimensional space, and since $\log(t)$ is an $i$-dimensional vector, there
must be a constant $d$ depending on the function $S$ such that for each $l$ we have
$$-d\log(C_i) < \log(t_l) < d\log(C_i) \eqno (2.21a)$$
Equation $(4.11a)$ is equivalent to 
$$C_i^{-e} < t_l < C_i^e \eqno (2.21b)$$
In particular, the variables $t_l$ are bounded away from 0.  
Next, continuing to focus on the $i > 0$ case, we examine how the $x$ to $(s,t)$ coordinate change 
affects $W_{ij}$ in the first $n-i$ variables.
It turns out that the relevant inequalities are those provided by Lemma 2.0.
This lemma says that if $x \in W_{ij}$, $w$ is in the vertex set $v(S)$ of $N(S)$ and on the
face $F_{ij}$, and $w' \in v(S)$ but $w' \notin F_{ij}$, then we have
$$x^{w' - w} < (C_{i+1})^{-\delta}$$
Writing in $y$ coordinates, this becomes
$$y^{L(w' - w)}< (C_{i+1})^{-\delta} \eqno (2.22a)$$
We would like to encapsulate the condition that $x \in (0,\eta)^n$ through an equation analogous to
$(2.22a)$. Shrinking $\eta$ if necessary, we can assume that for each $m$, $x_m = x^{{\bf e}_m} < 
(C_{i+1})^{-\delta}$, and we express this in $y$ coordinates as
$$y^{L({\bf e}_m)}< (C_{i+1})^{-\delta} \eqno (2.22b)$$
Writing $L = (L_1,L_2)$ and $y = (s,t)$ like before,  equations $(2.22)$ become
$$s^{L_1(w' - w)} < (C_{i+1})^{-\delta} t^{L_2(w - w')} \eqno (2.23a)$$
$$s^{L_1({\bf e}_m)} < (C_{i+1})^{-\delta} t^{L_2(-{\bf e}_m)} \eqno (2.23b)$$
Equation $(4.11b)$ says that each component of $t$ 
is between $C_i^{-e}$ and $C_i^e$. So there is a constant $d'$ depending only $N(S)$ such that in
$(2.23)$ one has
$$C_i^{-e'} < t^{L_2(w - w')} < C_i^{e'} \eqno (2.24a)$$
$$C_i^{-e'} < t^{L_2(-{\bf e}_m)} < C_i^{e'} \eqno (2.24b)$$
So as long as $A_2$ from the beginning of section 3 is sufficiently large, equations $(2.23)$ give
$$s^{L_1(w' - w)} < 1 \eqno (2.25a)$$
$$s^{L_1({\bf e}_m)} < 1 \eqno (2.25b)$$
Summarizing, if $x \in W_{ij}$, then the corresponding $(s,t)$ in $y$ coordinates satisfy $(2.19b)$ and
$(2.25a)-(2.25b)$.
We now use in a similar fashion the other inequalities of Lemma 2.0. Namely, $x \in (0,\eta)^n$ is in 
$W_{ij}$ if $(2.19a)$ holds and $x$ satisfies the following for all $w \in v(S) \cap F_{ij}$, 
$w' \in v(S) \cap (F_{ij})^c$
$$x^{w'} < C_n^{-1}x^w \eqno (2.26a)$$
Analogous to above, we incorporate the condition $x \in (0,\eta)^n$ by stipulating that $\eta < (C_n)^
{-1}$ and write
$$x^{{\bf e}_m} < C_n^{-1} \eqno (2.26b)$$
Analogous to $(2.23)$, these can be written as 
$$s^{L_1(w' - w)} < (C_{n})^{-1} t^{L_2(w - w')} \eqno (2.27a)$$
$$s^{L_1({\bf e}_m)} < (C_{n})^{-1} t^{L_2(-{\bf e}_m)} \eqno (2.27b)$$
Again using $(2.24)$, there is some $\mu$ such that equations $(2.27)$ hold whenever for all $w' - w$
and all ${\bf e}_m$ we have
$$s^{L_1(w' - w)} < \mu  \eqno (2.28a)$$ 
$$s^{L_1({\bf e}_m)} < \mu  \eqno (2.28b)$$
Hence if a point $(s,t)$ is such that $s$ satisfies $(2.28a)-(2.28b)$ and $t$ satisfies $(2.19a)$, then
the corresponding $x$ is in $W_{ij}$. Putting $(2.25)$ and $(2.28)$ together, 
let $Y_{ij}$ denote the set $W_{ij}$ in the $y$ coordinates. Let $u_1$, $u_2$,... be an enumeration
of the set of all $L_1(w' - w)$ for vertices $w \in F_{ij}$ and $w' \notin F_{ij}$, as well as the 
distinct $L_1({\bf e}_m)$. 
We define the sets $E_1$ and $E_2$ by
$$E_1 = \{ s: 0 < s^{u_l} < \mu  \hbox { for all } l \} \times D_{ij} \eqno (2.29a)$$
$$E_2 = \{ s: 0 < s^{u_l} < 1  \hbox { for all } l \} \times D_{ij} \eqno (2.29b)$$
Then by $(2.25)$ and $(2.28)$ we have 
$$E_1 \subset Y_{ij} \subset E_2 \eqno (2.29c)$$
It is worth pointing out that none of the $u_l$ are zero: If some $\bar{w}_l - \bar{w}_0$ were 
zero this would imply that they came from a $w \in F_{ij}$ and a $w' \notin F_{ij}$ such that
$w' - w$ is a function of only the $t$-variables. This would mean that $w' - w$ is tangent to
$F_{ij}$, which can never happen when $w \in F_{ij}$ and $w' \notin F_{ij}$. If some 
$L_1({\bf e}_m)$ were zero, that would imply ${\bf e}_m$
is a function of the $t$ variables only, meaning that ${\bf e}_m$ is tangent to $F_{ij}$. Since 
$F_{ij}$ is a bounded face, this cannot happen either.

Equations $(2.29a)-(2.29c)$ are for $i > 0$, and there are analogous equations when $i = 0$.
Fortunately, these require less effort to deduce; a coordinate change is not required.
There is a single vertex $v$ on a given $F_{0j}$.
Lemma 2.0 tells us that if $\mu$ is sufficiently small, if we define 
$$F_1 = \{ x \in (0,\eta)^n: x^{v'} < \mu x^v \hbox { for all } v' \in v(S) - \{v\}\}$$
$$F_2 = \{ x \in (0,\eta)^n: x^{v'}< x^v \hbox { for all } v' \in v(S) - \{v\}\}$$
Then we have $F_1 \subset W_{0j} \subset F_2$. To combine this with the $i > 0$ case, we
rename the $x$ variables $s$ and define $Y_{0j} = 
W_{0j}$. Let $\{u_l\}_{l > 0} $ be an enumeration of the $v' - v$ for
$v' \in v(S) - \{v\}$ as well as the unit coordinate vectors ${\bf e}_m$. When $i = 0$ define 
$$E_1 = \{ s:0 < s^{u_l} < \mu  \hbox { for all } l > 0\}$$
$$E_2 = \{ s:0 < s^{u_l} <  1 \hbox { for all } l > 0\}\eqno(2.30)$$
Then, shrinking $\mu$ to less than $\eta$ if necessary, like above we have $E_1 \subset Y_{0j} \subset 
E_2$.

In the remainder of this section, we consider the $i > 0$ and $i = 0$ cases together. We still have some
work to do. Namely, we would like to replace the sets $\{ s: 0 < s^{u_l} < \mu  \hbox { for all } l \}$
or $\{ s: 0 < s^{u_l} < 1 \hbox { for all } l \}$ by cubes. To this end, 
we will divide up $Y_{ij}$ in the $s$ variables into finitely many pieces. A 
coordinate change in the $s$ variables will be performed on each piece taking it to a set which is a 
positive curved quadrant. This is done as follows. For $i > 0$ let $E_1'$ and $E_2'$ be defined by
$$ E_1' = \{ s :0 < s^{u_l} < \mu \hbox { for all } l > 0\}$$
$$ E_2' = \{ s :0 < s^{u_l} < 1 \hbox { for all } l > 0\}$$
When $i = 0$, let $E_1' = E_1$ and $E_2' = E_2$.
Writing  $S = (S_1,..,S_{n-i}) = (\log(s_1),..,\log(s_{n-i}))$, in the $S$ coordinates $E_2'$ becomes
the set $E_2^S$ given by
$$E_2^S = \{ S: S \cdot u_l < 0  \hbox { for all } l\} $$
The set of $S$ satisfying $(2.30)$ is the intersection of several hyperplanes passing through 
the origin. We subdivide $E_2^S$ via the $n-i$ hyperplanes $S_m = 0$, resulting in (at most)
$2^{n-i}$ pieces which we call $E_2^{S,1}$, $E_2^{S,2}$,... We focus our attention on the one for which
all $S_m > 0$, which we assume is $E_2^{S,1}$. The intersection of $E_2^{S,1}$ with the hyperplane
$\sum_m S_m = 1$ is a polyhedron, which we can triangulate into finitely simplices $\{Q_p\}$ whose 
vertices all have rational coordinates. By taking the convex hull of these $Q_p$'s with the origin,
one obtains a triangulation of $E_2^{S,1}$ into unbounded $n$-dimensional
"simplices" which we denote by $\{R_p\}$. Each $R_p$ has $n$ unbounded faces of dimension $n-1$
 containing the origin. The equation for a given face can be written as $S \cdot q^{p,l} = 0$, where
each $q^{p,l}$ has rational coordinates, so that
$$R_p = \{S : S \cdot q^{p,l} < 0 \hbox { for all } 1 \leq l \leq n - i \} \eqno (2.31)$$
Hence $\cup R_p = E_2^{S,1}$. The other $E_2^{S,m}$ can be similarly subdivided. We 
combine all simplices from all the $E_2^{S,m}$ into one list $\{R_p\}$. 
Note each $R_p$ on the combined list satisfies $(2.31)$. Furthermore, the $R_p$ are disjoint and
 up to a set of measure zero $E_2^S = \cup_p R_p$. Converting back now into $s$
coordinates, for $i > 0$ we define 
$$Y_{ijp} = \{(s,t) \in Y_{ij}: \log(s) \in R_p\} = \{(s,t) \in Y_{ij}:0 < s^{q^{p,l}} < 1 \hbox 
{ for all } 1 \leq l \leq n - i \}\eqno (2.32a)$$
When $i = 0$ we let
$$Y_{0jp} = \{s \in Y_{0j}: \log(s) \in R_p\} = \{s \in Y_{0j}:0 < s^{q^{p,l}} < 1 \hbox 
{ for all } 1 \leq l \leq n \}\eqno (2.32b)$$
Then the $Y_{ijp}$ are disjoint and up to a set of measure zero we have
$$\cup_p Y_{ijp} = Y_{ij} \subset  E_2\eqno (2.33)$$
On each $Y_{ijp}$ we shift from $y = (s,t)$  coordinates 
(or $y = s$ coordinates if $i = 0$) to $z = (\sigma,t)$ coordinates (or $z = \sigma$ coordinates if
$i = 0$), where $\sigma$ is  defined by
$$\sigma_l = s^{q^{p,l}} \hbox { for } l \leq n - i\eqno (2.34)$$
In the new coordinates, $Y_{ijp}$ becomes a set $Z_{ijp}$ where 
$$ Z_{ijp} \subset (0,1)^{n-i} \times D_{ij}\,\,\,\,\,\,\,\,\,\,(i > 0) \eqno (2.35a)$$
$$ Z_{ijp} \subset (0,1)^n\,\,\,\,\,\,\,\,\,\,(i = 0) \eqno (2.35b)$$
Let $W_{ijp}$ denote the set $Z_{ijp}$ in the original 
$x$ coordinates. So the $W_{ijp}$ are disjoint open sets and up to a set of measure zero 
$\cup_p W_{ijp} = W_{ij}$.

\noindent {\bf Lemma 2.4.} If $i > 0$, write $z = (\sigma, t)$, where $\sigma$ denotes the first $n - i$
 components and $t$ the last $i$ components. For any vector $w$, we denote by $(w',w'')$ the vector 
such that the
monomial $x^w$ transforms to $\sigma^{w'}t^{w''}$ in the $z$ coordinates. In the case where $i = 0$,
write $z = \sigma$ and say that $x^w$ transforms into $\sigma^{w'}$.

\noindent {\bf a}) If $w$ is either a unit coordinate vector ${\bf e}_l$, or of the form $v' - v$ for 
$v$ a vertex of $S$ in $F_{ij}$ and $v'$ a vertex of $S$ not in $F_{ij}$, then each component of 
$w'$ is nonnegative, with at least one component positive.

\noindent {\bf b}) If each component of $w$ is nonnegative, then so is each component of $w'$ and $w''$.

\noindent {\bf c}) There exists some $\mu' > 0$ such that for all $i,j$, and $p$
$$(0,\mu')^{n-i} \times D_{ij} \subset Z_{ijp} \subset (0,1)^{n-i} \times D_{ij} \,\,\,\,\,\,\,\,\,\,(i > 0)\eqno (2.36a)$$
$$(0,\mu')^n \subset Z_{ijp} \subset (0,1)^n\,\,\,\,\,\,\,\,\,\,(i = 0)\eqno (2.36b)$$
In particular, when $i > 0$, for fixed $t$ the cross-section of $Z_{ijp}$ is 
a positive curved quadrant.

\noindent {\bf Proof.} We assume that $i > 0$; the $i = 0$ case is done exactly the same way. If 
$w$ is of one the forms of part a), then the monomial $x^w$ in the $x$ 
coordinates becomes a monomial of the form $s^{u_m}t^a$ in the $y$ coordinates, where the $u_m$ 
are as before.
Since $Y_{ijp} \subset E_2$, where $E_2$ is as in $(2.29)$ or $(2.30)$, whenever each $s^{q^{p,l}} < 1$ for each $l$
we have $s^{u_m} < 1$ for each $m$. Thus if we write $s^{u_m} = \prod_l (s^{q^{p,l}})^{\alpha_l}$, each
$\alpha_l$ must be nonnegative; otherwise we could fix any $s^{q^{p,l}}$ for which $\alpha_l$ is 
nonnegative, and let the remaining $s^{q^{p,l}}$ go to zero, eventually forcing $s^{u_m} = \prod_l
(s^{q^{p,l}})^{\alpha_l}$ to be greater than 1. This means that the $\alpha_l$ are nonnegative. If
they were all zero, this would mean $u_m = 0$ which cannot happen by the discussion after $(2.29c)$.
So at least one $\alpha_l$ is positive.
Since $s^{u_m}t^{v}$ transforms into $\sigma^{\alpha_l}t^v$ in the $z$ coordinates, we have part 
a) of this lemma.

Next, we saw that any $x_l$ transforms into some $s^{a_l}t^{b_l}$ in the $y$ coordinates, where each
component of $a_l$ and $b_l$ is nonnegative . When transforming from $x$ to $z$ coordinates, by part
a) $x_l$ transforms into some  $\sigma^{a_l'}t^{b_l}$ with $a_l'$ having nonnegative components. 
Hence part b) holds for the $x_l$. Therefore it holds for any $x^w$ with each
component of $w$ nonnegative.

Moving to part c), the right-hand sides follow from $(2.35)$. As for the left hand sides, from the 
expression $s^{u_m} = 
\prod_l (s^{q^{p,l}})^{\alpha_l}$ with nonnegative $\alpha_l$, there is a $\mu' > 0$ such that each 
$s^{u_m} < \mu$ whenever $s^{q^{p,l}} < \mu'$ for all $l$. So if $s^{q^{p,l}} < \mu'$ for 
each $l$ and $t \in D_{ij}$, then $(s,t) \in E_1$. By $(2.29c)$, we conclude that whenever 
$s^{q^{p,l}} < \mu'$ for all $l$ and if $t \in D_{ij}$, then $y = (s,t)$ is in $Y_{ijp}$. In the
$z$ coordinates this becomes the left hand inequality of $(2.36a)$ for $i > 0$. When $i = 0$, the 
same argument holds; 
whenever $s^{q^{p,l}} < \mu'$ for each $l$ then $s \in E_1$ and $(2.36b)$ follows. Thus we are done 
with the proof of Lemma 2.4.

\noindent We can now give the proof of Theorem 2.2.

\noindent {\bf Proof of Theorem 2.2.} Parts a) and b) follow from part c) of Lemma 2.4 except the
statement that $D_{ij} \subset (C_i^{-e}, C_i^e)^i$ which is a consequence of $(2.21b)$ and the fact
that the $y$ to $z$ coordinate changes do not affect the $t$ variables. Moving on to part c),
the discussion prior to $(2.19)$ showed that for $v_1$ and $v_2$ on $F_{ij}$, the $x$ to $y$ coordinate
change takes $x^{v_2 - v_1}$ to a function of the $t$ variables only. Since the coordinate change from $y$
to $z$ variables do not affect the $t$ variables, the  $x$ to $z$ coordinate
change takes $x^{v_2 - v_1}$ to a function of the $t$ variables only as well, giving that $v_1' = 
v_2'$ as required.

Next, if $v$ is a vertex of $N(S)$ on $F_{ij}$ and $v_l$ is a vertex of $N(S)$ not on $F_{ij}$ 
or is of the form $v + {\bf e}_m$ for some $m$, then by Lemma 2.4a) $(v_l')_k \geq (v')_k$ for
all $k$ with at least one component strictly positive. Any $w \in N(S)$ satisfies
$w - v = \sum c_l (v_l - v) + \sum d_m {\bf e}_m$ for some nonnegative $c_l$ and $d_m$, where
$v_l$ are vertices of $N(S)$. As long as $w \notin F_{ij}$, there is either
going to be some positive $c_l$ for $v_l \notin F_{ij}$, or some positive $d_m$. Hence in this
situation some $(w' - v')_k > 0$. This completes the proof of Theorem 2.2.
 
Suppose $x = f(z) = (z^{m_1},...,z^{m_n})$, where $m_i = (m_{i1},...,m_{in})$ are 
vectors such that ${\rm det}(m_{ij})$ is nonzero. Then a direct calculation reveals that the Jacobian 
determinant of $f(z)$ is given by
$${\rm det}(m_{ij}) (z^{\sum_i m_i - (1,1,...,1,1)}) \eqno (2.37)$$
If in addition all the $m_{ij}$ are nonnegative, we can find a $g(z)$ of the form $g(z) = (z_1^{k_1},...
z_n^{k_n})$, $k_j > 0$, such that $f \circ g (z)$ has constant determinant. To see this, one uses the 
chain rule in conjunction with $(2.37)$. One gets that the determinant of $f \circ g(z)$ is given by
$$(\prod_l k_l){\rm det}(m_{ij}) \prod_j z_j^{k_j \sum_i m_{ij} - 1}$$
Hence by setting $k_j = {1 \over \sum_i m_{ij}}$, one obtains that $f \circ g (z)$ has constant 
determinant. (The invertibility of $(m_{ij})$ insures that none of these sums are zero). Note that in
Theorem 2.2, if one replaces $\beta_{ijp}(z)$ by such a $\beta_{ijp} \circ g(z)$, the conclusions
of the theorem continue to hold. Hence in the rest of this paper, without losing generality we assume 
that for all $i$, $j$, and $p$, the Jacobian determinant of $\beta_{ijp}(z)$ is constant. One advantage
of doing this is that integrals transform simply under $\beta_{ijp}(z)$ this way. Another is illustrated
by the following lemma.

\noindent {\bf Lemma 2.5.} Suppose $z = h(x) = (x^{b_1},...,x^{b_n})$, where $b_i =  
(b_{i1},...,b_{in})$ is such that the determinant of $B = (b_{ij})$ is nonzero and the 
Jacobian determinant of $h(x)$ is constant. Let $\beta_i$ denote the hyperplane through the origin 
spanned by
the vectors $b_j$ for $j \neq i$. Then a monomial $x^{\alpha}$ transforms into the monomial $z^{\tilde
{\alpha}}$ in $z$ coordinates, where the $i$th component ${\tilde{\alpha}}_i$ is given by any 
component of the intersection of the hyperplane $\beta_i + \alpha$ with the line $\{(t,t,....,t,t):
t \in \R\}$

\noindent {\bf Proof.} We use the notation $B_{j,v}$ to denote the
matrix obtained by replacing the $j$th row of $B$ by the vector $v$. The hyperplane $\beta_j + \alpha$ has 
equation ${\rm det}(B_{j,x}) = {\rm det}(B_{j, \alpha})$, so a component of the intersection of this plane with the
line $\{(t,t,....,t,t): t \in \R\}$ is given by
$${{\rm det}(B_{j, \alpha}) \over {\rm det}(B_{j, (1,1,...,1,1)})}\eqno (2.38)$$
Next we examine how a monomial $x^{\alpha}$ transforms under the $x$ to $z$ coordinate change. To
understand this, we work in logarithmic coordinates. Writing $X = (\log(x_1),...,
\log(x_n))$ and $Z = (\log(z_1),..., \log(z_n))$, one has that $Z = BX$ or $X = B^{-1}Z$ where $X$ 
and
$Z$ are viewed as $n$ by 1 column matrices. The function $\log(x^{\alpha})$ becomes $\alpha^T X = 
\alpha^T B^{-1} Z$. Thus in the $z$ coordinates, $x^{\alpha}$ becomes $z^{\tilde{\alpha}}$, where
$\tilde{\alpha} = (B^T)^{-1}\alpha$.  By Cramer's rule, 
$(B^T)^{-1} \alpha = {1 \over {\rm det}(B)}({\rm det}(B_{1,\alpha}),...,{\rm det}(B_{n,\alpha}))$. 
Comparing with $(2.38)$,
to prove this lemma we must show that ${\rm det}(B_{j,(1,1,...,1,1)}) = {\rm det}(B)$ for all $j$. 

To accomplish
this, we use the fact that the Jacobian determinant of $h$ is a constant function. By $(2.37)$, this
means  we have $\sum_i b_i = (1,1,...,1,1)$. In matrix form, this can be written as 
$$(1,1,...,1,1)B  = (1,1,...,1,1) \eqno (2.39)$$
Writing ${\bf 1} = (1,1,...,1,1)$, taking transposes of $(2.39)$ gives
$$B^T {\bf 1} = {\bf 1}$$
Equivalently, 
$${\bf 1} = (B^T)^{-1}{\bf 1}$$
By Cramer's rule this means that for all $j$ we have
$$1 = {{\rm det}(B_{j,(1,1,...,1,1)}) \over {\rm det}(B)}$$
This is what we need to show and we are done.

Lemma 2.6 will interpret Lemma 2.5 in the setting of Theorem 2.2. To this end, let $W_{ijp}$ be one of the
open sets of Theorem 2.2 and $\beta_{ijp}:Z_{ijp} \rightarrow W_{ijp}$ be the associated map. We 
write $z = (\beta_{ijp})^{-1}(x) = (z^{b_1},...,z^{b_n})$. Here $b_i = (b_{i1},...,b_{in})$ where 
$(b_{ij})$ is an invertible matrix of rational numbers which can be negative.

Let $v$ be a vertex of $N(S)$ on a face $F_{ij}$. Write $z = (\sigma,t)$,
where $\sigma$ are the first $n-i$ coordinates and $t$ are the last $i$ coordinates. The monomial
$x^v$ transforms into some $\sigma^{v'}t^{v''}$ in the $z$ coordinates in accordance 
with Lemma 2.5. The $t^{v''}$ factor is of little interest; by Theorem 2.2 the $t$ coordinates
are bounded above and below away from zero and thus so is $t^{v''}$. The vector $v'$ on the other hand
is very important for the purposes of this paper, and Lemma 2.6 gives the relevant properties:

\noindent {\bf Lemma 2.6.} Let $v' = (v'_1,...,v'_{n-i})$ be as above. Let $d$ be the Newton distance
of $S$, and let $C(S)$
be the face of $N(S)$ (possibly unbounded) such that the line $\{(t,t,....,t,t): t \in \R\}$ intersects
$C(S)$ in its interior. Let $k$ be the dimension of $C(S)$, where $k = 0$ if the line intersects 
$N(S)$ at a vertex. Then the following hold.

\noindent {\bf a)} Each $v_m'$ satisfies $0 \leq v'_m \leq d$ \parskip = 3pt.

\noindent {\bf b)} At most $n - k$ of the $v'_m$ are equal to $d$. 

\noindent {\bf c)} If $n - k$ of the $v'_m$ are equal to $d$, then the face $F_{ij}$ is a subset of 
$C(S)$.

\noindent {\bf d)} If $F_{ij} = C(S)$, then all $n - k$ of the $v'_m$ are equal to $d$.
\parskip = 12pt

\noindent {\bf Proof.} Let $1 \leq m \leq n - i$. As in the proof of Lemma 2.5, a monomial 
$z^{\alpha}$ becomes $x^{B^T \alpha}$ in the $x$ coordinates. The image of the hyperplane 
$\{\alpha: \alpha_m = 0\}$ under
$ B^T$ is the span of the $b_l$ for $l \neq m$, denoted by $\beta_m$ in Lemma 2.5. Hence the
image of $\{\alpha: \alpha_m = v'_m\}$ is the hyperplane $\beta_m + v$. Suppose $w \in N(S)$. Let $x^w$ 
transform into $\sigma^{w'}t^{w''}$ in the $z$ coordinates. By Theorem 2.2c), 
$\sigma^{w'}$ has at least as high a power of $\sigma_m$ appearing as does $\sigma^{v'}$. In other 
words $w'_m \geq v'_m$. Translating back into the $x$ coordinates, any such $w$
must
be on a single side of the hyperplane $\beta_m + v$; we conclude that $\beta_m + v$ is a 
separating hyperplane for $N(S)$. Furthermore, by Theorem 2.2 c), for $v, w \in F_{ij}$ one has that
$w'_m = v'_m$. Translating this into the $x$ coordinates, we have that this
separating hyperplane in fact contains the face $F_{ij}$. 

Since $\beta_m + v$ is a separating hyperplane for $N(S)$, it cannot intersect the line $\{(t,t,...
.,t,t): t \in \R\}$ in the interior of $N(S)$. Thus the intersection point is some $(t,t,...,t,t)$ 
with $t \leq d$, with $t = d$ only if $(d,d,...,d,d) \in \beta_m + v$. Hence by Lemma 2.5, $v'_{m}
\leq d$. By Lemma 2.4 b) we also have $v'_{m} \geq
0$, so we conclude that $0 \leq v'_{m} \leq d$ for all $1 \leq m \leq n - i$, giving a).

We now analyze how many of the $v'_{m}$ can actually be equal to $d$. Let $p$ denote the 
number of $v'_{m}$ that are equal to $d$. By the above discussion, if $m$ satisfies $v'_{m} = d$, then
$\beta_m + v$ must contain $F_{ij}$ as well as the point $(d,d,...,d,d)$. Since any separating hyperplane
for $N(S)$ containing $(d,d,...,d,d)$ must also contain all of $C(S)$, we have that such a 
$\beta_m + v$
in fact contains $span(C(S), F_{ij})$. Hence the intersection of all $p$ of these $\beta_m + v$ contains
$span(C(S), F_{ij})$. We conclude that
$$n - p = dim(\cap_{\{m:v'_m = d\}} (\beta_m + v)) \geq dim(span(C(S), F_{ij}) \geq dim(C(S)) = k \eqno (2.40)$$
We conclude that $p \leq n - k$, giving b). Furthermore, if $p = n - k$, all the
inequalities in $(2.40)$ must be equalities. In particular, $dim(span(C(S), F_{ij}) =  dim(C(S))$. The
only way this can happen is if $F_{ij} \subset C(S)$, giving c). Lastly, suppose $F_{ij} = C(S)$. Then
since each hyperplane $\beta_m + v$ for $1 \leq m \leq n - i = n - k$ contains $F_{ij}$, each such 
hyperplane also contains $(d,d,...,d,d)$. Hence by Lemma 2.5, each $v'_m = d$ and we have part d). This 
concludes the proof.

\noindent {\bf 3. Proofs of lower bounds of Theorems 1.2 - 1.4}

\noindent We start with this elementary lemma, which we will make repeated use of. 

\noindent {\bf Lemma 3.1.} Suppose $m_1,...,m_n$ are nonnegative numbers not all zero. Let $M = \max_i
m_i$, and let $l$ denote the number of $m_i$ equal to $M$. Then if $|E|$ denotes Lebesgue measure,
we have the following for all $0 < \delta < 1$, where $C$ and $C'$ are constants depending on the 
$m_i$. 

\noindent {\bf a}) 
$$C |\ln \delta|^{l-1}\delta^{1 \over M} < |\{x \in (0,1)^n: x_1^{m_1}...x_n^{m_n} < \delta\}| < C' |\ln \delta|^{l-1}
\delta^{1 \over M} \eqno (3.1)$$
\noindent {\bf b)} If $M < 1$, then 
$$C \delta  < \int_{\{x \in (0,1)^n: {\delta \over x_1^{m_1}...x_n^{m_n}} < 1\}} {\delta \over
x_1^{m_1}...x_n^{m_n}}\,dx < C' \delta $$
\noindent {\bf c)} If $M = 1$, then 
$$C |\ln \delta|^l \delta  < \int_{\{x \in (0,1)^n:  {\delta \over x_1^{m_1}...x_n^{m_n}} < 1\}} {\delta \over x_1^{m_1}...x_n^{m_n}}\,dx
< C' |\ln \delta|^l \delta $$
\noindent {\bf d)} If $M > 1$, then 
$$\int_{\{x \in (0,1)^n: {\delta \over x_1^{m_1}...x_n^{m_n}} < 1\}} {\delta \over x_1^{m_1}...x_n^{m_n}}\,dx < 
C' |\{x \in (0,1)^n: x_1^{m_1}...x_n^{m_n} < \delta\}| \eqno (3.2)$$ 

\noindent {\bf Proof.} We first deal with parts b) and c). Note that when each $m_i \leq 1$, we have
$$(\delta^{1 \over n},1)^n \subset \{x \in (0,1)^n: {\delta \over x_1^{m_1}...x_n^{m_n}} < 1\} \subset 
(\delta,1)^n$$
Thus 
$$\int_{\{x \in (0,1)^n: (\delta^{1 \over n},1)^n\}} {\delta \over x_1^{m_1}...x_n^{m_n}}\,dx < \int_{\{x \in (0,1)^n: {\delta 
\over  x_1^{m_1}...x_n^{m_n}} < 1\}} {\delta \over x_1^{m_1}...x_n^{m_n}}\,dx  $$
$$< \int_{(\delta,1)^n} {\delta \over x_1^{m_1}...x_n^{m_n}}\,dx \eqno (3.3)$$
One can integrate the left and right hand sides of $(3.3)$ directly and get parts b) and c). Moving on
to a), we proceed by induction on $n$. When $n = 1$ it is immediate, so assume $n > 1$ and the result is
known for $n - 1$. Without losing generality, we may assume that $m_n = M$. We regard $|\{x \in (0,1)^n: 
x_1^{m_1}... x_n^{m_n} < \delta\}|$ as the integral of the characteristic function of 
$\{x \in (0,1)^n: x_1^{m_1}... x_n^{m_n} < \delta\}$, integrating with respect to $x_n$ first. We have
$$|\{x \in (0,1)^n: x_1^{m_1}... x_n^{m_n}< \delta\}| = \int_{(0,1)^{n-1}} \min (1, {\delta^{1/M} \over x_1^{m_1/M}....
x_{n-1}^{m_{n-1}/M}})\,dx \eqno (3.4)$$
Break $(3.4)$ into 2 parts, depending on whether or not $x_1^{m_1}....x_{n-1}^{m_{n-1}} < \delta$.
The portion where $x_1^{m_1}....x_{n-1}^{m_{n-1}} < \delta$ gives a contribution of 
$|\{x \in (0,1)^{n-1}: x_1^{m_1}....x^{m_{n-1}} < \delta\}|$, which by induction hypothesis will always
be smaller by at least a factor of $C |\ln \delta|$ than the left and right hand sides of $(3.1)$. As
for the the portion where $x_1^{m_1}....x_{n-1}^{m_{n-1}} > \delta$, one obtains the integral
$$\int_{\{x \in (0,1)^{n-1}:  {\delta^{1/M} \over x_1^{m_1/M}...x_{n-1}^{m_{n-1}/M}} < 1\}} 
({\delta^{1/M} \over x_1^{m_1/M}...x_{n-1}^{m_{n-1}/M}})\,dx \eqno (3.5)$$
Since $M \geq m_i$ for all $i < n$, one can estimate $(3.5)$ using parts b) or c) of this lemma. Since
exactly $l - 1$ of ${m_1/M},...,{m_{n-1}/M}$ are equal to 1, if $l > 1$ part c) says that $(3.5)
\sim \delta^{1/M} |\ln \delta|^{l-1}$ as needed, while if $l = 1$ part b) says that $(3.5) \sim
\delta^{1/M}$ as needed. This completes the proof of a).

Moving on to d), we again may assume that $m_n = M$ and perform the $x_n$ integration first. We have
$$\int_{\{x \in (0,1)^n: {\delta \over x_1^{m_1}...x_n^{m_n}} < 1\}} {\delta \over x_1^{m_1}..
.x_n^{m_n}}$$
$$= \int_{\{x \in (0,1)^{n-1}: {\delta \over x_1^{m_1}...x_{n-1}^{m_{n-1}}} < 1\}}
(\int_{1 > x_n > {\delta^{1/M}  \over x_1^{m_1/M}...x_{n-1}^{m_n/M}}} {\delta \over x_1^{m_1}...x_n^{m_n}}
\, dx_n)dx_1...dx_{n-1}$$
Since $m_n > 1$, this is bounded by
$$C \int_{\{x \in (0,1)^{n-1}: {\delta \over x_1^{m_1}...x_{n-1}^{m_{n-1}}} < 1\}}
(\int_{{2\delta^{1/M}  \over x_1^{m_1/M}...x_{n-1}^{m_{n-1}/M}} > x_n > {\delta^{1/M}  
\over x_1^{m_1/M}...x_{n-1}^{m_{n-1}/M}}} {\delta \over x_1^{m_1}...x_n^{m_n}}\, dx_n)$$
$$dx_1...dx_{n-1} \eqno (3.6)$$
The integrand is bounded above by a constant, so this is at most
$$C|\{x \in (0,1)^{n-1}\times (0,2): {x_1^{m_1}...x_{n}^{m_n}} < 2^M\delta \}| \eqno (3.7)$$
Rescaling in the $x_n$ variable and using part a) gives us part d) and we are done. 

We now start the proofs of the lower bounds of Theorems 1.2-1.4. 
Note that the lower bounds of Theorem 1.3 are contained in those of Theorem 1.2, so it 
suffices to prove the lower bounds of Theorem 1.2 to prove both.

\noindent {\bf Proof of Theorem 1.2a)} Let $R(x) = \sum_{v \in v(S)} x^v$, where as earlier in this 
paper $v(S)$ denotes the set of vertices of $N(S)$. Note that $R(x)$ and $S(x)$ have the same Newton
polyhedron. By the corollary to Lemma 2.1, there is a constant $C$ such that $|S(x)| \leq C|R(x)|$ 
for all $x \in (0,\infty)^n$. Hence it suffices to show Theorem 1.2a) for $|R|$ in place of $|S|$. 

\noindent {\bf Case 1)} The face $C(S)$ is compact.

\noindent Let $F_{kj}$ denote $C(S)$, and let $W_{kjp}$ the corresponding sets from Theorem
2.2. We have $I_{|R|,\phi}(\epsilon) = \int_{|R| < \epsilon} \phi(x)\,dx$. Note that it suffices to show 
that each
$\int_{\{x \in W_{kjp}: |R| < \epsilon\}} \phi(x)\,dx > C |\ln\epsilon|^{n - k - 1} 
\epsilon^{1 \over d}$ for some constant $C$. The $x$ to
$z$ coordinate change has constant Jacobian determinant by the discussion above Lemma
2.5, so if $\Phi$ denotes $\phi \circ \beta_{kjp}$, where $\beta_{kjp}$ is as in Theorem 2.2 we have 
$$\int_{\{x \in W_{kjp}: |R| < \epsilon\}} \phi(x)\,dx = c\int_{\{z\in Z_{kjp}: |R \circ \beta_{kjp}| < 
\epsilon\}} \Phi(z)\,dz$$
Since $\phi(0) > 0$, $\Phi(0) > 0$ as well, so for some $\delta, \xi > 0$ we have  
$$\int_{\{z\in Z_{kjp}: |R \circ \beta_{kjp}| < \epsilon\}} \Phi(z)\,dz > \delta |\{z\in Z_{kjp} \cap
(0,\xi)^n: |R \circ \beta_{kjp}(z)| < \epsilon\}| $$
By part a) of Theorem 2.2, we have $(0,\mu')^n \subset Z_{kjp}$ for some $\mu' > 0$. Hence for $\rho
= \min(\mu',\xi)$ we have
$$\int_{\{x \in W_{kjp}: |R| < \epsilon\}} \phi(x)\,dx  > C|\{ z \in (0,\rho)^n: |R \circ \beta_{kjp}(z)| < 
\epsilon\}|$$
Writing $z = (\sigma,t)$ as in Theorem 2.2, each function $x^v$ for $v \in v(S)$ transforms into some 
function $\sigma^{v'}t^{v''}$ in the $z$ coordinates where the components of $v'$ and $v''$ are all nonnegative.
By part c) of Theorem 2.2, each component of $v'$ is minimized for $v \in F_{kj} = 
C(S)$, and by part b), each $t_i$ is bounded above and below away from zero. Hence if we fix some
$V \in F_{kj}$, for $z \in Z_{kjp}$ we have
$$|R \circ \beta_{kjp}(z)| = |\sum_v \sigma^{v'}t^{v''}| < C \sigma^V \eqno (3.8)$$
We conclude that 
$$|\{ z \in (0,\rho)^n: |R \circ \beta_{kjp}(z)| < \epsilon\}| > C''|\{z \in (0,\rho)^n : 
C|\sigma^V| < \epsilon\}| \eqno (3.9)$$
By part d) of Theorem 2.6, each component of $V$ is just equal to $d$. So by Lemma 3.1 a) (scaled), we have
$$|\{z \in (0,\rho)^n : C|\sigma^V| < \epsilon\}| > C'' |\ln\epsilon|^{n - k - 1} \epsilon^{1 \over d}$$
This gives the desired lower bounds and we are done in case 1.

\noindent {\bf Case 2)} The face $C(S)$ is unbounded. Let $V = \sum_{l=1}^n a_l x_l = c$ denote a
separating hyperplane for $N(S)$ such that $V \cap N(S) = C(S)$. Note that each $a_l$ is
nonnegative. Since $C(S)$ is unbounded, at least one $a_l = 0$. Without loss of generality, we may
let $q < n$ such that $a_l > 0$ for $1 \leq l \leq q$ and $a_l = 0$ for $l > q$. Correspondingly
write $x = (x',x'')$, where $x' \in \R^q$ and $x'' \in \R^{n-q}$. Let $P$ denote the 
projection onto the first $q$ coordinates. Define $\bar{R}(x') = \sum_{v \in P(v(S))} (x')^v$. From
first principles one can verify that
$$P(N(R)) = N(\bar{R})$$
Using that $V$ is a separating hyperplane for $N(R)$ it is also straightforward to verify that $P(V)$ is a
separating hyperplane for $P(N(R)) = N(\bar{R})$ with $N(\bar{R}) \cap P(V) = P(C(S))$.  But the 
equation for 
$P(V)$ is given by $\sum_{l=1}^q a_l x_l = c$ and each $a_l > 0$ for $l \leq q$. Thus $P(C(S))$ is
a compact face of $N(\bar{R})$. Furthermore, since the directions ${\bf e}_l$ for $l > q$ are all 
parallel to $C(S)$ and $(d,d,...,d,d) (n$ times) is in the interior of $C(S)$, $(d,d,...,d,d) 
(p $ times) is in the interior of $P(C(S))$. For the same reasons, the codimension of $P(C(S))$ in
$\R^q$ is the same as the 
codimension of $C(S)$ in $\R^n$, namely $n - k$. Hence we may apply Case 1 to $\bar{R}(x')$ and get
the lower bounds of Theorem 1.2, for $\bar{R}(x')$ in place of $S(x)$. 

For a given $v \in v(S)$, we write $v = (v',v'')$ where $v'$ denotes the first $q$ components and 
$v''$ the last $n - q$ components. We can write
$$I_{|R|,\phi}(\epsilon) = \int_{\R^{n-q}}(\int_{\{x' \in \R^{q}: |R(x',x'')| < \epsilon\}}
\phi(x',x'')\,dx')\,dx'' \eqno (3.10)$$
Since $\phi(0) > 0$, there are $\delta, \xi > 0$ such that $(3.10)$ is greater than
$$\delta \int_{(0,\xi)^{n-q}}|\{x' \in (0,\xi)^{q}: |R(x',x'')| < \epsilon\}|\,dx'' \eqno (3.11)$$
For fixed $x'$, one has $R(x',x'') = \sum_{v \in v(S)} (x')^{v'}(x'')^{v''}$, so consequently
$$|R(x',x'')| < C \sum_{v' \in P(v(S))} (x')^{v'}  = C\bar{R}(x') \eqno (3.12)$$
Hence by $(3.10)$ and $(3.11)$ we have
$$I_{|R|,\phi}(\epsilon) > \delta \xi^{n-q}|\{x' \in (0,\xi)^{q}: |\bar{R}(x')| < {\epsilon \over 
C}\}|\eqno (3.13)$$
As indicated above, case 1) of this lemma applies to $\bar{R}(x')$, which has the same values of 
$d$ and $k$ that $R(x)$ (and $S(x)$) do. Choosing an appropriate $\phi$ we get 
$$|\{x' \in (0,\xi)^{q}: |\bar{R}(x')| < {\epsilon \over C}\}| > C' |\ln\epsilon|^{n - k - 1} 
\epsilon^{1 \over d} \eqno (3.14)$$
Combining $(3.13)$ and $(3.14)$ gives the desired result and we are done.

To prepare for the proof of the lower bounds of Theorem 1.4, we consider the setting of 
Theorem 2.2, focusing on a specific $i$, $j$, and $p$. For now assume that $i > 0$. Note that
$\beta_{ijp}$ is defined on all of $[0,\infty)^n$, not just $Z_{ijp}$. Furthermore, by Lemma 2.4a), for 
any
$t \in [0,\infty)^i$, $\beta_{ijp}(0,t) = 0$. Hence $S \circ \beta_{ijp}$ is defined on a neighborhood 
of $\{0\} \times [0,\infty)^i$ in $[0,\infty)^n$.
Write $S(x) = \sum s_{\alpha}x^{\alpha}$ like before. By Theorem 2.2c), there 
is a single $\omega$ such that if $\alpha \in F_{ij}$, $x^{\alpha}$ transforms in the $z$ 
coordinates into $\sigma^{\omega}t^{\alpha''}$ for some $\alpha''$ that depends on $\alpha$ Hence $S_{F_{ij}}(x) = 
\sum_{\alpha \in F_{ij}} s_{\alpha}x^{\alpha}$ transforms into $\sigma^{\omega}P(t)$, where $P(t)$
is a polynomial in $t^{{1\over N}}$ for some $N$. Any of our conditions on $S_{F_{ij}}(x)$ 
translates into a corresponding
condition on $P(t)$. On $Z_{ijp}$ we may write
$$S \circ \beta_{ijp}(z) = \sigma^{\omega}P(t) + \sum_{\alpha \notin F_{ij}} s_{\alpha}
\sigma^{\alpha'}t^{\alpha''} \eqno (3.15)$$
Equation $(3.15)$ assumed that $i >0$, but the $i= 0$ case can be incorporated by letting $t = 1$
and letting $P(1)$ be the appropriate coefficient. Using Theorem 2.2c) again, for a given $\alpha$ 
in the sum $(3.15)$ each $\alpha'_k \geq \omega_k$ with at least one
inequality strict. Since $\sum s_{\alpha}x^{\alpha}$ is a convergent Taylor series, we have
$|s_{\alpha}| < CR^{|\alpha|}$ for some $C$ and $R$. Because of Lemma 2.4b), $|\alpha|$ and 
$|\alpha'| + |\alpha''|$ are within a constant factor of one another. Hence we have an estimate
$|s_{\alpha}| < C'(R')^{|\alpha'| + |\alpha''|}$ and therefore for some $N$ the sum in $(3.15)$ 
represents a Taylor series in $\sigma_k^{1 \over N}$ and $t_k^{1 \over N}$ convergent near the
origin, not just on $Z_{ijp}$. Consequently, for some real-analytic functions $r_k(\sigma,t)$ of $\sigma_k^{1 \over N}$ 
and $t_k^{1 \over N}$ we can rewrite $(3.15)$ as 
$$S \circ \beta_{ijp}(z) = \sigma^{\omega}[P(t) + \sum_{k=1}^{n-i} (\sigma_k)^{1 \over N}
r_k(\sigma,t)]\eqno (3.16)$$
Equation $(3.16)$ is valid near the origin. But it is also valid on a neighborhood
of $\{0\} \times [0,\infty)^i$ in $[0,\infty)^n$. (If $i = 0$, we take 
$[0,\infty)^i$ to mean $\{1\}$). To see this, note that for any $\beta \leq N\omega$, we have
$\partial_{\sigma}^{\beta} (S \circ \beta_{ijp}((\sigma_1^N,...,\sigma_{n-i}^N,t) - \sigma^{N\omega}P(t))$ is 
zero on a set $\{0\} \times U$ where $(3.16)$ is known to hold. Hence by real-analyticity it must 
be true on all of $\{0\} \times [0,\infty)^i$. This 
implies that $(3.16)$ makes sense on a neighborhood of $\{0\} \times [0,\infty)^i$ in $[0,\infty)^n$.

We now proceed to the proof of Theorem 1.4. Assume $C(S)$ is compact face of codimension $k$, and 
there is some $x' \in (\R - \{0\})^n$ such that the growth index of $|S_{C(S)}|$ at $x'$ is 
$a  \leq {1 \over d}$ with multiplicity $q \geq 0$. Without loss of generality we may assume that
$x' \in (0,\infty)^n$. Let $F_{kj} = C(S)$, and let
$W_{kjp}$ and $Z_{kjp}$ be any of the sets of Theorem 2.2 corresponding to this face. In the $z$
coordinates, $S_{C(S)}(x)$ becomes $\sigma^{\omega}P(t)$. Under this coordinate change, $x'$ becomes
some $z' = (\sigma',t')$ where $P(t)$ has growth index $a$ at $t'$ with multiplicity $q$.
Since the coordinate change has constant determinant, if $\beta_{kjp}$ as in Theorem 2.2 and
$\Phi$ denotes $\phi \circ \beta_{kjp}$, then 
$$I_{|S|,\phi}(\epsilon) = \int_{|S| < \epsilon}\phi(x)\,dx \geq  \int_{\{x \in (\R^+)^n: |S(x)| 
< \epsilon\}}\phi(x)\,dx $$
$$=  c \int_{\{z \in (\R^+)^n:|S\circ \beta_{kjp}(z) < \epsilon|} \Phi(z)\,dz \eqno (3.17)$$
Since $\beta_{kjp}(0,t) = 0$ for all $t$ by Lemma 2.4a), $\Phi(0,t) = \phi(0) > 0$ for all $t$. Thus 
we may let
$U$ be a neighborhood of $(0,t')$ in $(\R^+)^n$ such that $\Phi(\sigma,t) > {\phi(0) \over 2}$ 
on $U$. We then have
$$I_{|S|,\phi}(\epsilon) > {\phi(0) \over 2} |\{z \in U: |S\circ \beta_{kjp}(z)| < \epsilon\}|
\eqno (3.18)$$
Hence it suffices to find a lower bound for $|\{z \in U: |S\circ \beta_{kjp}(z)|< \epsilon\}|$
We will do this by finding a lower bound for 
$$|\{z = (\sigma,t) \in (\epsilon^{{1\over d(n - k)}},\epsilon^{{1\over d(n - k) }+ \mu})^{n-k} \times U' : |S\circ \beta_{kjp}(z)|
<  \epsilon\}| \eqno (3.19)$$
Here $U'$ is a neighborhood of $t'$, and
$\mu$ is a sufficiently small positive number to be determined. We may assume $\epsilon$ is small
enough that $(3.16)$ holds on the set in $(3.19)$. Using $(3.16)$, we rewrite $(3.19)$ as
$$|\{(\sigma,t) \in (\epsilon^{{1\over d(n - k) }},\epsilon^{{1\over d(n - k)} + \mu})^{n-k} \times U
': |P(t) + \sum_{l=1}^{n-k} (\sigma_l)^{1 \over N}r_l(\sigma,t)| < {\epsilon \over \sigma^{\omega}}\}
| \eqno (3.20)$$
By Lemma 2.6d), each $\omega_l = d$, so $(3.20)$ is just
$$|\{(\sigma,t) \in (\epsilon^{{1\over d(n - k) }},\epsilon^{{1\over d(n - k)} + \mu})^{n-k} \times U
': |P(t) + \sum_{l=1}^{n-k} (\sigma_l)^{1 \over N}r_l(\sigma,t)| < {\epsilon \over \sigma_1^d....\sigma_{n-k}^d}\}|
\eqno (3.21)$$
When $\sigma \in (\epsilon^{{1\over d(n - k)}},\epsilon^{{1\over d(n - k)} + \mu})^{n-k}$, one has that 
${\epsilon \over \sigma_1^d....\sigma_{n-k}^d} >  \epsilon^{d\mu(n-k)}$. On the other hand,
$\sum_{l=1}^{n-k} (\sigma_l)^{1 \over N}r_l(\sigma,t) < C\epsilon^{{1\over dN(n - k)} + 
{\mu \over dN}}$. Thus if $\mu$ were chosen appropriately small, then for small enough $\epsilon$,
if $\sigma \in (\epsilon^{{1\over d(n - k)}},\epsilon^{{1\over d(n - k)} + \mu})^{n-k}$ one has
$$ |\sum_{l=1}^{n-k} (\sigma_l)^{1 \over N}r_l(\sigma,t)| < {1 \over 2}{\epsilon \over \sigma_1^d....
\sigma_{n-k}^d} \eqno (3.22)$$
Consequently, for such $\epsilon$, $(3.20)$ is bounded below by
$$|\{(\sigma,t) \in (\epsilon^{{1\over d(n - k) }},\epsilon^{{1\over d(n - k)} + \mu})^{n-k} 
\times U': |P(t)|  < {\epsilon \over 2\sigma_1^d....\sigma_{n-k}^d}\}| $$
$$= \int_{(\epsilon^{{1\over d(n - k) }},\epsilon^{{1\over d(n - k)} + \mu})^{n-k}}
|\{t \in U': P(t)  < {\epsilon \over 2\sigma_1^d....\sigma_{n-k}^d}\}|\,d\sigma \eqno (3.23)$$
By virtue of the facts that $t' \in U'$ and $P(t)$ has growth index $a$ at $t'$ with 
multiplicity $q$, the integrand in $(3.23)$ is bounded below by 
$C(\ln|{\epsilon \over \sigma_1^d....\sigma_{n-k}^d}|)^q ({\epsilon \over \sigma_1^d....
\sigma_{n-k}^d})^a$. Hence $(3.23)$ is bounded below by
$$C\int_{(\epsilon^{{1\over d(n - k) }},\epsilon^{{1\over d(n - k)} + \mu})^{n-k}}
(\ln|{\epsilon \over \sigma_1^d....\sigma_{n-k}^d}|)^q ({\epsilon \over \sigma_1^d.... \sigma_{n-k}^d})
^a\,d\sigma \eqno (3.24)$$
Scaling each of the $\sigma$ variables by $\epsilon^{{1\over d(n - k) }}$, $(3.24)$ becomes
$$C\epsilon^{1 \over d}\int_{(1,\epsilon^{-\mu})^{n-k}}
(\ln(\sigma_1....\sigma_{n-k}))^q ({1 \over \sigma_1^{da}.... \sigma_{n-k}^{da}})
\,d\sigma \eqno (3.25)$$
We now evaluate $(3.25)$ on a case by case basis. If $a = {1 \over d}$, one can do a term by term 
expansion of the logarithm in the integrand of 
$$C\epsilon^{1 \over d}\int_{(1,\epsilon^{-\mu})^{n-k}}
(\ln(\sigma_1) + .... + \ln(\sigma_{n-k}))^q ({1 \over \sigma_1^ .... \sigma_{n-k}})\,d\sigma \eqno 
(3.26)$$
Integrating $(3.26)$ term by term becomes immediate, and results in a lower bound of 
$$C |\ln\epsilon|^{q + n - k}\epsilon^{1 \over d}$$
This is the lower bound of Theorem 1.4b). On the other hand if $a < {1 \over d}$, we may choose 
$f$ with $a < f < {1 \over d}$,
and we have
$$(\ln(\sigma_1....\sigma_{n-k}))^q ({1 \over \sigma_1^{da}.... \sigma_{n-k}^{da}}) > 
C{1 \over \sigma_1^{df}.... \sigma_{n-k}^{df}}$$ 
Hence it suffices to find lower bounds for 
$$\epsilon^{1 \over d}\int_{(1,\epsilon^{-\mu})^{n-k}} {1 \over \sigma_1^{df}.... 
\sigma_{n-k}^{df}}\,d\sigma \eqno (3.27)$$
This is easily integrated directly to give a lower bound
$$C\epsilon^{{1 \over d} - (n-k)\mu(1 - df)}$$
Setting $a' = {1 \over d} - (n-k)\mu(1 - df)$ gives Theorem 1.4a) and we are done.

\noindent {\bf 4. Proofs of upper bounds of Theorems 1.2 and 1.3.} Recall that 
$$I_{|S|, \phi}(\epsilon) =  \int_{\{x : |S(x)| < \epsilon\}} \phi(x) \,dx$$
We will bound $\int_{\{x \in (\R^+)^n : |S(x)| < \epsilon\}} \phi(x) \,dx$ as the other octants are
entirely analogous. We may assume that $\phi$ is supported in $(-\eta,\eta)^n$ where $\eta$ is as in
the constructions of section 2. Since $\phi$ is bounded, it suffices to bound a given
$$|\{x  \in (0,\eta)^n: |S(x)| < \epsilon\}| = \sum_{ijp}|\{x  \in W_{ijp} : |S(x)| < \epsilon\}|
$$
Clearly it is enough to bound each term separately. Since for each $i,j,$ and $p$ the $x$ to $z$
coordinate change has constant Jacobian, it suffices to bound
$$|\{z  \in Z_{ijp} : |S \circ \beta_{ijp} (z)| < \epsilon\}| \eqno (4.0)$$
So our task is to bound $(4.0)$ by the appropriate right hand side of Theorems 1.2 and 1.3. 
We now fix some $i$,$j$, and $p$. Let $a$ denote the 
maximum order of any zero of $S_{F}(x)$ on $(\R - \{0\})^n$, for any compact face $F$ of $N(S)$. 
In the notation of $(3.15)$, 
this implies that the order of any zero of $P(t)$ on $(\R - \{0\})^i$ is at most $a$. By well
known methods (see [S] Ch 8 sec 2.2), this means for any $t \in (\R - \{0\})^i$, there is some
directional derivative $\partial_w$ and some $0 \leq a' \leq a$ such that $\partial_w^{a'}P$
is nonzero. (If $i = 0$ we take $a' = 0$). Note that by Theorem 2.2b) if $(\sigma,t) \in Z_{ijp}$ then $t \in (C_i^{-e},C_i^e)^n$.
By continuity and compactness, we can let $\{E_l\}$ be a finite collection of cubes covering
$[C_i^{-e},C_i^e]^n$, 
$w_l$ be directions, $a_l$ be nonnegative integers, and $\delta_0 > 0$ a constant such that
on $E_l$ 
$$|\partial_{w_l}^{a_l}P(t)| > \delta_0 \eqno (4.1)$$
We next examine the effect of taking such directional derivatives on the sum in $(3.15)$. Using the
fact that $|\alpha''| < C|\alpha|$ for some $C$, taking any
$t$ directional derivative of order at most $a$ on this sum leads to a term bounded by 
$$C\sum_{\alpha \notin F_{ij}} |s_{\alpha}||\alpha|^a \sigma^{\alpha'}t^{\alpha''}(\min_m t_m)^{-a} 
\eqno (4.2)$$
We may assume that the $E_l$ are small enough so that $t_m > {1 \over 2} C_i^{-e}$ for each $m$
on each $E_l$. Hence $(4.2)$ is bounded by 
$$C'C_i^{ae} \sum_{\alpha \notin F_{ij}} |s_{\alpha}||\alpha|^a \sigma^{\alpha'}t^{\alpha''} 
\eqno (4.3)$$
By Lemma 2.1, if $(\sigma,t) \in Z_{ijp}$, then for some $V \in F_{ij}$ $(4.3)$ is bounded by 
$$C'C_i^{ae}C_{i+1}^{-\delta}x^V = C'C_i^{ae}C_{i+1}^{-\delta}\sigma^{V'}t^{V''} = 
C'C_i^{ae}C_{i+1}^{-\delta}\sigma^{\omega}t^{V''}\eqno (4.4)$$
Here $\omega$ is as in $(3.15)$.
We can assume $|t_l| < 2 C_i^{e}$ for each $l$, so for some $e'$ equation $(4.4)$ is bounded by
$$C'C_i^{ae'}C_{i+1}^{-\delta}\sigma^{\omega}\eqno (4.5)$$
We can assume $C_{i+1}$ was chosen small enough so that $C'C_i^{ae'}C_{i+1}^{-\delta} < {\delta_0 
\over 2}$; shrinking $C_{i+1}$ has no effect on any of the coordinate changes for the $i$-dimensional
faces, or on the constant $C'C_i^{ae'}$ in $(4.5)$. Hence we can assume that $(4.5)$ is 
bounded by
$${\delta_0 \over 2} \sigma^{\omega}\eqno (4.6)$$
Combining $(4.1)$ and $(4.6)$ in $(3.15)$, we conclude that for $(\sigma,t) \in Z_{ijp}$ with 
$t \in E_l$ one has 
$$| \partial_{w_l}^{a_l} (S \circ \beta_{ijp}(z))| > {\delta_0 \over 2}
\sigma^{\omega}\eqno (4.7)$$
We now prove the appropriate bounds $(4.0)$. Note that it suffices to bound each
$$|\{z = (\sigma,t)  \in Z_{ijp} : t \in E_l, |S \circ \beta_{ijp} (z)| < \epsilon\}|\,\,\,\,\,\,\,(i > 0)\eqno (4.8a)$$
$$|\{z = \sigma  \in Z_{0jp} : |S \circ \beta_{0jp} (z)| < \epsilon\}|\,\,\,\,\,\,\,(i = 0)\eqno (4.8b)$$
To do this, we separate into cases $a_l = 0$ and $a_l > 0$. For $a_l = 0$, by $(4.7)$, equation $(4.8a)$
or $(4.8b)$ is at most 
$$ C|\{\sigma \in (0,1)^{n-i} : \sigma^w < {2 \over \delta_0} \epsilon\}|\eqno (4.9)$$
By Lemma 2.6a), each component of $\omega$ is at most the Newton distance $d$, and the number of
times $d$ may appear in $\omega$ is at most the codimension $n - k$ of the face called $C(S)$. 
Hence by Theorem 3.1a), we have that $(4.9)$ is at most
$$C'|\{\sigma \in (0,1)^{n-i}: \sigma^{\omega} < \epsilon\}| < C'|\ln\epsilon|^{n - k - 1}
\epsilon^{{1 \over d}} \eqno (4.10)$$
This term is no greater than any of the right hand sides in Theorem 1.2, so we do not have
to worry about it any further. We now move to the case when $a_l > 0$. Here we use Van der Corput's 
lemma in the $w_l$ direction and then integrate the result. Since the $Z_{ijp}$ are defined through
monomial inequalities, their cross-sections in the $w_l$ direction consist of boundedly many segments.
Applying the van der Corput lemma $(2.1)$ of [C1], we see 
that the $w_l$ cross section of $(4.8a)$ has measure at most $C({\epsilon \over \sigma^{\omega}})^
{{1 \over a_l}} =  {\epsilon^{{1 \over a_l}}\over \sigma^{\omega/a_l}}$. Here $\omega/a_l$ denotes the 
vector where each component of $\omega$ is divided by $a_l$. It also of course has measure at 
most $C$ since the $t$ variables are bounded. Hence  $(4.8)$ is bounded by
$$C \int_{(0,1)^{n-i}} \min(1, {\epsilon^{{1 \over a_l}}\over \sigma^{\omega/a_l}})\,d\sigma
\eqno (4.11)$$
It is natural to divide $(4.9)$ depending on whether or not ${\epsilon \over \sigma^{\omega}} < 1$.
We get that $(4.11)$ is bounded by
$$C|\{\sigma \in (0,1)^{n-i}: \sigma^{\omega} < \epsilon\}| + C \int_{{\epsilon \over \sigma^{\omega}} < 1}
{\epsilon^{{1 \over a_l}}\over \sigma^{\omega/a_l}}\,d\sigma $$
The left hand term is exactly $(4.10)$ and satisfies the desired bounds in all cases. Since each $a_l$ is at most the maximum order $a$ of any zero of any 
$S_{F}(x)$, the second term of $(4.10)$ is at most
$$C \int_{{\epsilon \over \sigma^{\omega}} < 1}{\epsilon^{{1 \over a}}\over \sigma^
{\omega/a}} = C \int_{{\epsilon^{1 \over a} \over \sigma^{\omega/a}} < 1}{\epsilon^{{1 \over a}}\over \sigma^
{\omega/a}}\,d\sigma \eqno (4.12)$$
To analyze $(4.12)$, we use the various parts of Lemma 3.1 to obtain the various upper 
bounds of Theorem 1.2. First suppose $a < d$. Then one or more components of $\omega/a$ may be 
greater than one. If this is in fact the case, Theorem 3.1d) says that $(4.12)$ is bounded by
the expression $(4.10)$, which is the needed bound of the second statement of Theorem 1.2b).
If all components of $\omega/a$ are at most 1, then by Theorem 3.1b) or c), $(4.12)$ is at most 
$C|\ln\epsilon|^{n-i}\epsilon^{1 \over a}$. Since $a < d$, this is better than the bound 
$C|\ln\epsilon|^{n- k - 1}\epsilon^{{1 \over d}}$ required by the second
statement of Theorem 1.2b). This completes the proof of Theorem 1.2 for $a < d$.

If $a = d$,
then each component of $\omega/a$ is at most 1, with at most $n - k$ equal to 1, so by Theorem 3.1c),
$(4.12)$ is at most $|\ln\epsilon|^{n-k}\epsilon^{{1 \over d}}$. This is the bound needed for the
first statement of Theorem 1.2b). By Lemma 2.6c), the only way $n - k$ components of $\omega/a$
could be equal to 1 is for $F_{ij}$ to be a subset of $C(S)$. If this is not the case, then Lemma
3.1c) says that $(4.12)$ is at most $C|\ln\epsilon|^{n- k - 1}\epsilon^{{1 \over d}}$. For a subface
of $C(S)$ with zeroes of order at most $b < a = d$, then as in the $a < d$ case $(4.11)$ is at most 
$C|\ln\epsilon|^{n- k - 1}\epsilon^{{1 \over d}}$.
Hence as long as $C(S)$ has no compact subface $F$ such that $S_F(x)$ has a zero of order $d$, one
gets the upper bound $C|\ln\epsilon|^{n- k - 1}\epsilon^{{1 \over d}}$ of the second statement of Theorem 1.2b). Thus we have proven Theorem 1.2 for 
$a = d$.

If $a > d$, then each component of $\omega/a$ is less than 1, so by 
Theorem 3.1b) $(4.12)$ is bounded by $C\epsilon^{{1 \over a}}$, the bound needed for Theorem 1.2 c)
and we are done.

We now move on to the proof of the upper bounds in Theorem 1.3. As in the proof for Theorem 1.2, it
suffices to prove upper bounds for 
$$\int_{\{x \in (0,\eta)^n : |S(x)| < \epsilon\}} \phi(x) \,dx = \sum_{ijp}\int_{\{x \in W_{ijp} : 
|S(x)| < \epsilon\}} \phi(x) \,dx$$
$$= \sum_{ijp}c_{ijp}\int_{\{z \in Z_{ijp} : |S\circ \beta_{ijp}(z)| < \epsilon\}} \Phi_{ijp}(z) 
\,dz\eqno (4.13)$$
Here $\Phi_{ijp}(z)$ denotes $\phi \circ \beta_{ijp}(z)$ and $c_{ijp}$ is the (constant) Jacobian 
determinant of the $x$ to $z$ coordinate change. Clearly, it suffices to prove upper bounds for
a given term of $(4.13)$. The proof of Theorem 1.2 carries through when $i < 2$ since the nondegeneracy
assumptions of Theorems 1.2 and 1.3 are the same for vertices and 1-dimensional edges and this is what
was used in the analysis of the $i < 2$ terms. Hence the estimates of Theorem 1.2 hold for those 
terms, which
imply the desired upper bounds in Theorem 1.3. So we assume that $i = 2$. Thus there are one $\sigma$
variable and two $t$ variables.

Let $D_{2j}$ be as in Theorem 2.2. Fix $t' \in cl(D_{2j})$. We may let
$U \times V$ be a neighborhood of $(0,t')$ in $[0,\infty)^3$ such that the expression 
$S \circ \beta_{ijp}(z) = \sigma^{\omega}[P(t) +  \sigma^{1 \over N}r(\sigma,t)
]$ of $(3.16)$ is valid on $U \times V$. 
Let $a$ denote the infimum over all compact faces $F$ of $N(S)$ and all $x \in (R - \{0\})^3$ 
of the growth index of $S_F$ at $x$. Since the $x$ to $z$ coordinate change transforms $S_{F_{2j}}(x)$ 
into $\sigma^{\omega}P(t)$, the infimum of the growth indices of $P(t)$ on $(R - \{0\})^2$ is at least
$a$. In particular, if we denote the growth index of $P(t)$ at $t = t'$ by $a(t')$, we have
$$a(t') \geq a \eqno (4.14)$$
In particular if $P(t') = 0$, then for a fixed $\mu > 0$ one has
$$a(t') > a - \mu \eqno (4.15)$$
So in this situation, if  $V$ is  sufficiently small, which we may assume, for any $\epsilon > 0$ 
we have
$$|\{t \in V: |P(t)| < \epsilon\}| < C\epsilon^{a(t') - \mu}$$
Furthermore, by a stability theorem of Karpushkin [K], if $U$ is sufficiently small, which we may 
also assume, when each $\sigma_k \geq 0$ we have
$$|\{t \in V: |P(t) + \sigma^{1 \over N}r(\sigma,t)| < \epsilon\}| < 
C\epsilon^{a(t') - \mu} \eqno (4.16)$$
(Technically Karpushkin's result applies to analytic functions of $\sigma$ not $\sigma^{1 \over
N}$, but a simple change of variables in $\sigma$ gives us what we need).
Using compactness, we may let $\{U_l \times V_l\}$ be a finite collection of $U \times V$ 
covering $\{0\} \times cl(D_{2j})$ such that for a given $l$ either $P(t)$ doesn't vanish on
$cl(V_l)$, or $P(t)$ has a zero on $V_l$ with $(4.16)$ holding for $\sigma \in U_l$. Since the
continuous $\beta_{2jp}$ takes $\{0\} \times 
[0,\infty)^2$ to the origin, and other points of $[0,\infty)^3$ to points other than the origin,
if the support of $\phi$ is sufficiently small, which we may assume, then the support of $\Phi = 
\phi \circ \beta_{2jp}$
is contained in the neighborhood $\cup_l(U_l \times V_l)$ of $\{0\} \times cl(D_{2j})$. Hence to
bound $(4.13)$ it suffices to bound each
$$\int_{\{(\sigma,t)  \in U_l \times V_l: |S\circ \beta_{2jp}(z)| < \epsilon\}}\Phi_{2jp}(z) \,dz$$
Since $\Phi_{2jp}(z)$ is bounded, this is at most
$$C|\{(\sigma,t) \in U_l \times V_l: |S \circ \beta_{2jp}(z)| < \epsilon\}| \eqno (4.17)$$
For the $U_l \times V_l$ for which $P(t)$ doesn't vanish on $cl(V_l)$, one is in the setting
of Theorem 1.2; namely $(4.7)$ holds with $w_l = 0$ and the analysis there leading to $(4.10)$
gives bounds as strong as all right-hand sides of Theorem 1.3. Hence we may restrict our attention
to $l$ for which $P(t)$ has a zero in $V_l$. In this case, $(4.17)$ is at most
$$C|\{(\sigma,t) \in U_l \times V_l: |P(t) + \sigma^{1 \over N}r
(\sigma,t)| < {\epsilon \over \sigma^{\omega}}\}|$$
$$= \int_{U_l}|\{t \in V_l: |P(t) + \sigma^{1 \over N}r(\sigma,t)| < 
{\epsilon \over \sigma^{\omega}}\}|\,d\sigma \eqno (4.18)$$
Let $a'$ be the minimum of all the $a(t')$ corresponding to the different $U_l \times V_l$. So
in particular $a' \geq a$, where $a$ is as in $(4.14)$.
By the above-mentioned stability result of Karpushkin, the integrand of $(4.18)$ is at most
$C{\epsilon^{a' - \mu} \over \sigma^{\omega(a' - \mu)}}$. It
is also uniformly bounded by the measure of $V_l$. Hence $(4.18)$ is at most
$$ C\int_{U_l}\min(1, {\epsilon^{a' - \mu} \over \sigma^{\omega(a' - \mu)}})\,d\sigma \eqno (4.19)$$
It is natural to break up the integral $(4.19)$ into two parts, depending on whether or not
$|{\epsilon \over \sigma^{\omega}}|$ is less than or greater than 1. One gets that $(4.19)$ is bounded by
$$C|\{\sigma \in U_l: \sigma^{\omega} < \epsilon\}| + 
C\int_{ \{\sigma \in (0,1): {\epsilon^{a' - \mu} \over \sigma^{\omega(a' - \mu)}} < 1\}} 
{\epsilon^{a' - \mu} \over \sigma^{\omega(a' - \mu)}}\,d\sigma \eqno (4.20)$$
By Lemma 2.6, $\omega \leq d$. Thus the first term of $(4.20)$ is bounded by $C \epsilon^
{{1 \over d}}$. This is bounded by all the right hand sides of Theorem 1.3, so we need
only consider the second term of $(4.20)$. 

Consider the situation where $a \leq {1 \over d}$.
Then since $a' \geq a$, this second term of $(4.20)$ is bounded by 
$$C\int_{ \{\sigma \in (0,1): {\epsilon^{a - \mu} \over \sigma^{\omega(a - \mu)}} < 1\}} 
{\epsilon^{a - \mu} \over \sigma^{\omega(a - \mu)}}\,d\sigma \eqno (4.21)$$
Since $\omega \leq d$, we have that $\omega(a - \mu) < 1$. Thus we can apply Lemma 3.1b) (or 
integrate  directly) to obtain that the right
term of $(3.20)$ is at most $C\epsilon^{a - \mu}$. We conclude that the growth index of 
$|S|$ is at least $a - \mu$. Since this is true for all sufficiently small $\mu$, we conclude
that the growth index of $|S|$ is at least $a$. This gives us the first statement of Theorem 
1.3b) as well as Theorem 1.3c), using that the multiplicity of this index is at most 2. 

Next, we move to the setting of the second statement of Theorem 1.3b); that is, where the growth index 
of each $|S_F(x)|$ is greater than ${1 \over d}$ at every point in $(\R - \{0\})^n$. In this case
$a'$ is the minimum of finitely many numbers greater than ${1 \over d}$, and therefore $a' > 
{1 \over d}$. Assume $\mu$ is small enough that 
$a' - \mu > {1 \over d}$. In this case it is possible that 
$\omega(a' - \mu) > 1$ regardless of what $\mu$ is. If this happens, we use Lemma 3.1d), and
obtain that the second term of $(4.20)$ is bounded by a constant multiple of the first term,
which as indicated above is bounded by all right-hand sides of Theorem 1.3. In the case that
each $\omega(a' - \mu) \leq 1$, we apply
Lemma 3.1b) or c) to obtain that the second term of $(4.20)$ is at most 
$C|\ln\epsilon|\epsilon^{a' - \mu}$. Since $a' - \mu > {1 \over d}$, this is a better
estimate than the right hand side of the first equation of Theorem 1.3b), and we are done.

\noindent {\bf 5. Proofs of Theorems 1.5 and 1.6.} 

\noindent We start with the proof of Theorem 1.5, where we are working in two dimensions. 

\noindent {\bf Lemma 5.1.} If $F$ is a 1-dimensional compact edge of $N(S)$ not intersecting the
critical line $y = x$ in its interior, then $S_F(x)$ cannot have 
any zeroes on $(R - \{0\})^2$ of order greater than the Newton distance $d$. 

\noindent {\bf Proof.} Without loss of generality we assume $F$ lies entirely on or below the line 
$y = x$. Denote by $cx^ay^b$ the term of $S_F(x,y)$ with highest power of $y$ appearing. The line 
containing $F$ is a separating line for $N(S)$, so it intersects $N(S)$ at some $(d',d')$ for 
$d' \leq d$. it has negative slope, so $b \leq d' \leq d$. Since ${\partial_y^b}S_F(x,y) = 
cb!x^a$, we have a partial derivative of $S_F(x,y)$ of order at most $d$ that doesn't vanish on
$(R - \{0\})^2$. This completes the proof.

We now can prove Theorem 1.5. If the critical line doesn't intersect $N(S)$ in the interior of a
compact edge, then by Lemma 5.1 we are in the setting of the second statement of Theorem 1.2b). So 
the growth index of $S$ is ${1 \over d}$ and its multiplicity is $1 - k$. Hence the conculsions of
Theorem 1.5 are satisfied. 

Suppose now the critical line does intersect $N(S)$ in the interior of a
compact edge $F$. If the associated $S_F(x)$ has zeroes of order less than $d$, then Lemma 5.1
implies we are once again in the setting of the second statement of Theorem 1.2b), and thus
Theorem 1.5 is again satisfied. If $S_F(x)$ has a zero of order $d$ but not greater, Theorem
1.4b) now says we have a growth index of ${1 \over d}$ but multiplicity $1$. In other words, the final
statement of Theorem 1.5b) is satisfied. If $S_F(x)$ has a zero of order greater than $d$, then
by Theorem 1.4a) the growth index of $S$ is less than ${1 \over d}$. Hence the last statement of
Theorem 1.5a) is verified, and we are done.

\noindent We now turn to the proof of Theorem 1.6. As in equation $1.4a$ we write
$$I_{S, \phi}(\epsilon) \sim \sum_{j = 0}^{\infty} \sum_{i=0}^{n-1} c_{ij}(\phi) \ln(\epsilon)^i 
\epsilon^{r_j} \eqno (5.1a)$$
Similarly, write 
$$I_{-S, \phi}(\epsilon) \sim \sum_{j = 0}^{\infty} \sum_{i=0}^{n-1} C_{ij}(\phi) \ln(\epsilon)^i 
\epsilon^{R_j} \eqno (5.1b)$$
We now no longer assume that $\phi$  has to be nonnegative. Recall that
$$J_{S, \phi}(\lambda) = \int_{\R^n} e^{i \lambda S(x)} \phi(x)\,dx \eqno (5.2)$$
Doing the integration of $(5.2)$ by first 
integrating over level sets $S = t$ and then with respect to $t$, one gets
$$\int_0^{\infty} {d I_{S, \phi}(t) \over d t} e^{i \lambda t}\gamma(t)\,dt + 
\int_0^{\infty} {dI_{-S, \phi}(t) \over d t} e^{-i \lambda t}\gamma(t)\,dt 
\eqno (5.3)$$
Here $\gamma(t)$ is a bump function equal to 1 on the range of $S$. One can differentiate $(5.1a)$
termwise, insert the result into $(5.3)$, and then integrate termwise (we refer to [G2] for details).
One obtains an expression
$$\sum_{j = 0}^{\infty} \sum_{i=0}^{n-1} c_{ij}'(\phi) \int_0^{\infty}\ln(t)^i t^{r_j-1}e^{i\lambda t}
\gamma(t)\,dt + 
\sum_{j = 0}^{\infty} \sum_{i=0}^{n-1} C_{ij}'(\phi)\int_0^{\infty}\ln(t)^i t^{R_j-1}e^{-i\lambda t} 
\gamma(t)\,dt \eqno (5.4)$$
It is well-known (see [F]) that for any $l > 0$, any real $\lambda$ one has 
$$\int_0^{\infty} e^{i \lambda t}\ln(t)^mt^{\alpha} \gamma(t)\,dt = {\partial^m \over \partial
\alpha^m} {\Gamma(\alpha+1) \over (-i\lambda)^{\alpha + 1}} + O(\lambda^{-l}) \eqno (5.5)$$
The dominant term of $(5.5)$ as $\lambda \rightarrow +\infty$ is given by
${\Gamma(\alpha+1) \ln(\lambda)^m \over (-i\lambda)^{\alpha + 1}}$
Next, note that the leading term of $(5.1a)$ or $(5.1b)$ will translate into the leading term 
of the asymptotic expansion for $(5.2)$ unless their corresponding terms cancel out in $(5.4)$.
The leading terms of $(5.1a)$ and $(5.1b)$ will be at most the 
term corresponding to the growth index of $|S|$. If there is any cancellation in $(5.4)$, then the 
result will be even faster decay for $J_{S,\phi}$. Hence the upper bounds of Theorem 1.2, 1.3, and
1.5 hold for $J_{S,\phi}$.

Suppose now $\phi(x)$ is a nonnegative function. It is not hard to check using $(5.5)$ that the
leading terms of the two series of $(5.4)$ are given by $c_{ij}(\phi)r_j{\Gamma(r_j) 
\ln(\lambda)^i \over (-i\lambda)^{r_j}}$ and
$C_{i'j'}(\phi)R_j{\Gamma(R_j) \ln(\lambda)^{i'} \over (i\lambda)^{R_{j'}}}$, where $c_{ij}\ln(t)^it^{r_j}$
and $C_{i'j'}\ln(t)^{i'}t^{R_j'}$ are the leading terms of $(5.1a)$ and $(5.1b)$. They can only cancel
out if $i = i'$ and $r_j = R_{j'}$. The numbers $c_{ij}$ and $C_{ij'}$ are then both 
positive since the integrals they come from are of nonnegative functions. Hence for there to be
cancellation, the ratio of $(-i\lambda)^{r_j}$ and $(i\lambda)^{r_j}$ must be a negative number.
For this to happen, $r_j$ must be an odd integer. We conclude that so long as the
growth index of $|S|$ is not an odd integer, the oscillatory index of $S$ is the same as 
this growth index. This implies that the results of Theorems 1.2-1.3 will hold for the oscillatory 
index. Furthermore, if $d > 1$ there will be no cancellation and therefore all of the statements 
analogous to Theorems 1.2-1.5 will hold 
for the oscillatory index. Similarly,
if $S$ does not take both positive and negative values in every neighborhood of the origin, 
then either $(5.1a)$ or $(5.1b)$ will be zero. Then there cannot be any cancellation; the growth
index of $S$ or $-S$ directly translates into the oscillatory index. Thus all of the statements analogous to 
Theorems 1.2-1.5 will hold for $J_{S,\phi}$. This completes the proof of Theorem 1.6.

\noindent {\bf 6. References.}

\noindent [AGV] V. Arnold, S Gusein-Zade, A Varchenko, {\it Singularities of differentiable maps
Volume II}, Birkhauser, Basel, 1988.\parskip = 3pt\baselineskip = 3pt

\noindent [BM] E. Bierstone, P. Milman, {\it Resolution of singularities in Denjoy-Carleman 
classes.} Selecta Math. (N.S.) {\bf 10} (2004), no. 1, 1-28. 

\noindent [C1] M. Christ, {\it Hilbert transforms along curves. I.
Nilpotent groups}, Annals of Mathematics (2) {\bf 122} (1985), no.3, 575-596.

\noindent [C2] M. Christ, {\it Convolution, curvature, and combinatorics, a
case study}, International Math. Research Notices {\bf 19} (1998) 1033-1048.

\noindent [DKo] J-P. Demailly, J. Kollar, {\it Semi-continuity of complex singularity exponents and 
Kahler-Einstein metrics on Fano orbifolds}, Ann. Sci. École Norm. Sup. (4) {\bf 34} (2001), no. 4, 
525-556. 

\noindent [GrSe] A. Greenleaf, A. Seeger, {\it Oscillatory and Fourier
integral operators with degenerate canonical relation}, Publicacions
Matematiques special issue: Proceedings of the El Escorial Conference 2000
 (2002), 93-141.

\noindent [F] M.V. Fedoryuk, {\it The saddle-point method}, Nauka, Moscow, 1977.

\noindent [IM] I. Ikromov, D. M\"uller, {\it On adapted coordinate systems}, to appear, Trans. Amer. 
Math. Soc.

\noindent [K] V. N. Karpushkin, {\it Uniform estimates for volumes}, Tr.
Math. Inst. Steklova {\bf 221} (1998), 225-231.

\noindent [G1] M. Greenblatt, {\it A Coordinate-dependent local resolution of singularities and 
applications},  J. Funct. Anal.  {\bf 255}  (2008), no. 8, 1957-1994.

\noindent [G2] M. Greenblatt, {\it Resolution of singularities, asymptotic expansions of oscillatory 
integrals, and related Phenomena}, submitted.

\noindent [G3] M. Greenblatt, {\it Newton polygons and local integrability
of negative powers of smooth functions in the plane}, Trans. Amer. Math. Soc. 
{\bf 358} (2006), no. 2, 657-670.

\noindent [G4] M. Greenblatt, {\it A direct resolution of singularities for functions of two variables 
with applications to analysis},  J. Anal. Math. {\bf 92}  (2004), 233-257.

\noindent [G5] M. Greenblatt, {\it Sharp $L\sp 2$ estimates for one-dimensional 
oscillatory integral operators with $C\sp \infty$ phase.} Amer. J. Math. 
{\bf 127} (2005), no. 3, 659-695.

\noindent [PS] D. H. Phong, E. M. Stein, {\it The Newton polyhedron and
oscillatory integral operators}, Acta Mathematica {\bf 179} (1997), 107-152.

\noindent [PSSt] D. H. Phong, E. M. Stein, J. Sturm, {\it On the growth and 
stability of real-analytic functions}, Amer. J. Math. {\bf 121} (1999), no. 3, 519-554.

\noindent [PSt1] D. H. Phong, J. Sturm, {\it Algebraic estimates, stability of local zeta functions, and 
uniform estimates for distribution functions},  Ann. of Math. (2)  {\bf 152}  (2000),  no. 1, 277-329.

\noindent [PSt2] D. H. Phong, J. Sturm, {\it On the algebraic constructibility of varieties of integrable 
rational functions on $C\sp n$}, Math. Ann. {\bf 323} (2002), no. 3, 453-484. 

\noindent [R] V. Rychkov, {\it Sharp $L^2$ bounds for oscillatory
integral operators with $C^\infty$ phases}, Math. Zeitschrift, {\bf 236}
(2001) 461-489.

\noindent [Se] A. Seeger, {\it Radon transforms and finite type conditions},
Journal of the American Mathematical Society {\bf 11} (1998) no.4, 869-897.

\noindent [S] E. Stein, {\it Harmonic analysis; real-variable methods,
orthogonality, and oscillatory integrals}, Princeton Mathematics Series Vol. 
43, Princeton University Press, Princeton, NJ, 1993.

\noindent [V] A. N. Varchenko, {\it Newton polyhedra and estimates of
oscillatory integrals}, Functional Anal. Appl. {\bf 18} (1976), no. 3, 
175-196.

\noindent [Va] V. Vassiliev, {\it The asymptotics of exponential integrals, Newton diagrams, and
classification of minima}, Functional Analysis and its Applications {\bf 11} (1977) 163-172. 
\parskip = 12pt \baselineskip = 12pt

\noindent Department of Mathematics \hfill \break
\noindent 244 Mathematics Building \hfill \break
\noindent University at Buffalo \hfill \break
\noindent Buffalo, NY 14260 \hfill \break\parskip = 3pt

\noindent Fields Institute \hfill \break
\noindent 222 College Street\hfill \break
\noindent Toronto, Ontario M5T 3J1 \hfill \break

\noindent email: greenbla@uic.edu
 \end